\documentclass{article}
\usepackage[cp1250]{inputenc}
\usepackage{amsfonts}
\usepackage{amsmath}
\usepackage{amssymb,latexsym}
\begin{document}

\newcommand{\be}{\begin{equation}}
\newcommand{\ee}{\end{equation}}

\large
\begin{center}
{\bf \Large H\"older continuity of the pluricomplex Green function
\\ and Markov brothers'  inequality \vskip 3mm }

\vskip 8mm

{\Large Miros\l aw Baran \ {\large and} \ Leokadia
Bialas-Ciez{\normalsize$^{^1}$}}
 \footnotetext[1]{corresponding author}

\vskip 5mm

Jagiellonian University, Faculty of Mathematics and Computer Science\\
Institute of Mathematics, 30-348 Krak\'ow, \L ojasiewicza 6, Poland\\
tel.: +48 12 664 66 34, \ fax: +48 12 664 66 75\\
{{\it e-mail:} Miroslaw.Baran@im.uj.edu.pl \ \ \ \ \ \ \ }\\
 {${\phantom{p}}$ \ \ \ \ \ \ \ \ Leokadia.Bialas-Ciez@im.uj.edu.pl}

\end{center}

\vskip 10mm \noindent {\bf Abstract.} Let $V_E$ be the
pluricomplex Green function associated to a compact subset $E$ of
$ \mathbb{C}^N$. The well known H\"older Continuity Property of
$E$ means that there exist constants $B>0,\gamma \in (0,1]$ such
that $V_E(z)\leq B\,\textrm{dist}(z,E)^\gamma$. The main result of
this paper says that this condition is equivalent to a Vladimir
Markov type inequality, i.e. $\|D^\alpha P\|_E\leq M^{|\alpha |}
(\deg P)^{m|\alpha |}\left(|\alpha |!\right)^{1-m}\|P\|_E$, where
$m,M>0$ are independent of the polynomial $P$ of $N$ variables. We
give some applications of this equivalence and we present its
generalization related to a notion of a fit majorant. Moreover, as
a consequence of the main result we obtain a criterion for the
H\"older Continuity Property in several complex variables of the
type of Siciak's L-regularity criterion.

\noindent

\vskip 10mm

\noindent {\it Keywords and phrases:}  \\ pluricomplex Green's
function, H\"older continuity property, Markov inequality

\vskip 6mm

\noindent {\it 2010 Mathematics Subject Classification:} 41A17,
32U35.

\vskip 2cm

{\bf  1. Introduction}

\vskip 3mm

Let $E$ be a compact set in $\mathbb{C}^N$. The pluricomplex
Green's function (with pole at infinity) of $E$ can be defined by
\[ V_E(z) := \sup\{ u(z)\: : \: u\in \mathcal{L}_N \ \ {\rm and} \ \
u\le 0 \ {\rm on} \ E \}, \ \ \ \ z\in\mathbb{C}^N,\] where
$\mathcal{L}_N$ is the Lelong class of all plurisubharmonic
functions in $\mathbb{C}^N$ of logarithmic growth at the infinity,
i.e.
\[\mathcal{L}_N\!:=\!\{u \in PSH(\mathbb{C}^N) \: :
\: u(z)-\log \|z\|_2 \le \mathcal{O}(1) \ \ {\rm as} \
\|z\|_2\rightarrow\infty\}\] (for background information, see
\cite{[12]}). Here $\|z\|_2$ stands for the Euclidean norm in $
\mathbb{K}^N,\ \mathbb{K}=\mathbb{C}$ or $ \mathbb{K}=\mathbb{R}$.
In the univariate case $V_E$ coincides with the Green's function
$g_{E}$ of the unbounded component of $\hat{\mathbb{C}}\setminus
E$ with logarithmic pole at infinity (as usual
$\hat{\mathbb{C}}=\mathbb{C}\cup\{\infty\}$).

If $V_E^*(z)$ is the standard upper regularization of $V_E$ then
it is well known (Siciak's theorem) that either $V_E^*\in
\mathcal{L}_N$ or $V_E^*\equiv +\infty$. It is also equivalent to
the fact that $E$ is a non-pluripolar or pluripolar set,
respectively. If we define the $L$-capacity of $E$ as
$C(E)=\liminf\limits_{||z||_2\longrightarrow\infty}
\frac{||z||_2}{\exp V_E^*(z)}$, then $E$ is a pluripolar set if
and only if $C(E)=0$.

A set $E$ is L-regular if $\lim\limits_{w\longrightarrow
z}V_E^*(w)=0$ for every $z\in E$. Siciak has proved that this is
equivalent to the continuity of $V_E$ in the whole space
$\mathbb{C}^N$. Therefore, L-regularity is one of the global
properties of $E$ and a crucial role is played here by the
continuity of $V_E$ near $E$.

Another global property of the set $E$ that depends only on the
behaviour of $V_E$ near $E$ is the H\"older continuity property
(HCP for short) of the pluricomplex Green's function $V_E$ (see
the result due to B\l ocki in \cite[Prop.3.5]{[21]}  or Prop.2.6
below). By Cauchy's inequality, one can prove that HCP implies the
{\it A.Markov inequality}, i.e. there exist constants $m\ge 1, \,
M>0$ such that for every polynomial $P$ of $N$ variables \be \|\,
{\rm grad}\, P\|_E \ \leq \ M \:({\rm deg}\,P)^{m} \: \|P\|_E.
\label{AMI} \ee If $E$ admits inequality (\ref{AMI}) then it is
said to be a {\it Markov set} and we write $E\in AMI(m, M$). To
reveal the importance of this property, we quote the following
result due to Ple\'sniak.

\vskip 3mm

\noindent \textbf{Theorem 1.1} (\cite[Th.3.3]{[18]}).  \textit{If
$E$ is a $\mathcal{C}^\infty$ determining compact subset of
$\mathbb{R}^N$ then the following statements are equivalent to
property (\ref{AMI}):}

(i) \textit{(Bernstein’s Theorem) If the distance of a function
$f: E \rightarrow \mathbb{R}$ from the space of polynomials of
degree at most $n$ forms a rapidly decreasing sequence (as
$n\rightarrow \infty$) then $f$ is a restriction to $E$ of a
function $\bar{f}\in \mathcal{C}^\infty(\mathbb{R}^N)$.}

(ii) \textit{The space $( \mathcal{C}^\infty(E),\tau_1)$ is
complete, where $\tau_1$ is the topology in
$\mathcal{C}^\infty(E)$ determined by the Jackson's seminorms.}

(iii) \textit{There exists a continuous linear operator $L:(
\mathcal{C}^\infty(E),\tau_1) \rightarrow
(\mathcal{C}^\infty(\mathbb{R}^N), \tau_0)$ such that $Lf_{|_E} =
f$ for each $f\in \mathcal{C}^\infty(E)$, where $\tau_0$ is the
natural topology in $\mathcal{C}^\infty(\mathbb{R}^N)$.}

(iv) \textit{There exist positive constants $C$ and $\mu$ such
that for every polynomial $P$ of degree at most $n$ \[
\|P\|_{E_{1/n^\mu}} \ \le \ C \: \|P\|_E,\] where $E_r=\{z\in
\mathbb{C}^N\: : \: {\rm dist}\, (z,E)\le r\}$.}

 \vskip 3mm

An exciting question is whether there exists a relationship
between the A.Markov inequality and the behaviour of the Green's
function near the considered set. It is known (see \cite{[5]})
that every Markov set $E\subset\mathbb{C} $ is not polar and $E$
is L-regular if $E\subset\mathbb{R}$ (\cite{[7]}). It seems that
A.Markov inequality (\ref{AMI}) implies H\"older continuity
property but a proof is an open problem mentioned e.g. in
\cite{[18]}. Actually, even the question about L-regularity of
Markov sets in the general case remains still open.

We shall make an attempt in the direction of solving this problem
by concentrating on a generalization of a inequality proved by
A.Markov's younger brother, V.Markov. He  discovered in 1892,
after a very detailed investigation, a precise but intricate
estimate for the $k$-th derivative of polynomials (see e.g.
\cite{[19]}): {\it for any polynomial $P$ of degree not greater
than $n$  \be \|P^{(k)} \|_{[-1,1]} \ \leq \ T_n^{(k)}(1) \ \|
P\|_{[-1,1]} \ = \ \frac{n^2[n^2-1]...[n^2-(k-1)^2]}{1\cdot 3
\cdot \ldots \cdot (2k-1)} \ \| P\|_{[-1,1]} \label{VM} \ee where
$T_n(x)=\cos (n\arccos x)$ is the $n$-th Chebyshev polynomial}
(for $k=1$ it was proved by A.Markov in 1889).

Inequality (\ref{VM}) inspired us to consider a new type of Markov
inequality (see Def.2.8 below). It turns out that this inequality
is equivalent to H\"older continuity property of the pluricomplex
Green's function. This is the main result of the paper (Th.2.12).

Although the definition of HCP is simple, its verification for
particular sets can be very complicated (see e.g. \cite{an,rr}).
The Carleson-Totik criterion (see \cite[Th.1.2, Th.1.7]{[10]})
merits mentioning here. It gives an equivalent condition for HCP
expressed in terms of capacities in a similar way to Wiener's
criterion for L-regularity. This criterion can be used for proving
HCP for a large family of sets. However, the Carleson-Totik
criterion holds only in the univariate complex case (or in
$\mathbb{R}^N$) and the equivalence is valid under certain
additional assumption on sets e.g. for sets satisfying an exterior
cone condition.  In this context, Th.2.12 of this paper provides a
useful tool for showing HCP especially when sets do not satisfy
the assumptions of the criterion mentioned above. We give some
examples of such an application of Th.2.12. Moreover, we prove a
rather surprising fact that it is sufficient to verify the
H\"older continuity property of $V_E$ only in $N$ canonical
directions (Cor.2.13). This allows us to show HCP for a large
class of sets.

\vskip 2mm

The paper is organized as follows. A statement of the main results
 is presented in Section~2. The next section contains proofs of these
results. In Section 4 we give a generalization of Th.2.12 related
to a notion of a fit majorant. The next section deals with compact
subsets of $\mathbb{R}^N \subset \mathbb{R}^N +i \mathbb{R}^N =
\mathbb{C}^N$, especially with convex bodies and also with UPC
sets, i.e. uniformly polynomially cuspidal sets. In the last
section we show some applications of Th.2.12 for disconnected
compact sets.

\vskip 5mm

{\bf  2. Notations and statement of the main results}

\vskip 3mm

The pluricomplex Green's function is closely related to
polynomials (see \cite{[20]} or \cite[Th.5.1.7]{[12]}) in view of
the formulas \be V_E (z) \  = \ \log  \Phi_E (z), \ \ \ \ z\in
\mathbb{C}^N , \label{V_E} \ee where $\Phi_E$ is the {\it Siciak
extremal function}, i.e. \[ \Phi_E(z) \ = \ \, \sup\left\{
|P(z)|^{1/n}/\|P\|_E^{1/n} \: : \: P\in \mathcal{P}
(\mathbb{C}^N), \ {\rm deg}\, P= n\ge 1, \ P|_E \not\equiv 0
\right\},\] $\mathcal{P} (\mathbb{K}^N)$ denotes the vector space
of polynomials of $N$ variables with coefficients in
$\mathbb{K}\in \{\mathbb{C},\mathbb{R}\}$ and $\| \cdot \|_E$ is
the maximum norm on $E$.

In order to investigate the behaviour of $V_E$ near $E$, we define
$$V_E^\bullet (z):=\sup\{V_E(x-w):\
x\in E, \ ||w||_2\leq ||z||_2\},\ z\in \mathbb{C}^N,$$ that is, a
radial modification of $V_E$. The definition and main properties
of $V_E^\bullet$ were presented by M.Baran, L.Bialas-Ciez, {\it
Comparison principles for compact sets in $ \mathbb{C}^N$ with HCP
and Markov properties} during the Conference on Several Complex
Variables on the occasion of Professor's J\'ozef Siciak's 80th
birthday, Krak\'ow, 4-8 July 2011.
 We set out (without proofs) the following
examples:

$\bullet$ \ if $E$ is a  unit ball in $\mathbb{C}^N$ (with respect
to a fixed complex norm) then $V_E^\bullet (z)=\log
(1+||z||_2/C(E))$,

$\bullet$ \ if $E$ is a convex symmetric body in $\mathbb{R}^N$
then $V_E^\bullet (z)=\log h(1+||z||_2/(2\,C(E)))$, where
$h(t)=t+\sqrt{t^2-1}$ for $t\geq 1$.

$\bullet$ \ if $E$ is a polar set then $V_E^\bullet (0)=0,\
V_E^\bullet|_{\mathbb{C}^N\setminus\{0\}}\equiv +\infty $.

For the non-polar sets we can obtain a very important fact which
is derived from Prop.1.4 in \cite{[14]} (cf.
\cite[Th.2.1c)]{[4]}): \vskip 3mm

\noindent \textbf{Proposition 2.1.} \textit{If $E$ is a
non-pluripolar compact subset of $\mathbb{C}^N$ and $$\rho_E (r) \
:=\ V_E^\bullet (z) \ \ \ \ \ for \ \ \ \|z\|_2=r,$$ then \
$t\longmapsto\rho_E (e^t)$ \ is an increasing convex function.}
\vskip 3mm

\noindent \textbf{Remark 2.2.} The function $\rho_E$ has the
following basic properties:

\vskip 1mm a) $\rho_{\lambda E}(r)=\rho_E(\lambda^{-1}r),\ \lambda
>0$,

b) $\rho_{E\times F}(r)=\max (\rho_E(r),\rho_F(r))$,

c) $\lim\limits_{r\rightarrow \infty}(\rho_E (r)-\log r)=-\log
C(E)$,

d) $\rho_E$ is increasing, continuous on $(0,+\infty)$ and
consequently, $0=\rho_E (0)\leq \lim\limits_{r\rightarrow
0+}\rho_E (r)$. Therefore, L-regularity is equivalent to the
equality $ \lim\limits_{r\rightarrow 0+}\rho_E (r)=0$.

\vskip 1mm

\noindent Indeed, equality a) can be checked by a standard
verification. Formula b) is a consequence of the well-known
product property of the pluripotential Green function. A behavior
of $\rho_E$ for $r$ near the infinity is related to the L-capacity
of $E$: like in the proof of Th.2.3 in \cite{[4]} we can show
equality c). Statement d) is deduced directly from Prop.2.1.

On can easily check that $\rho_{[-1,1]}(r)=\log h(1+r)$. In
general, it is rather difficult to calculate the exact values of
$\rho_E$. However, for some investigations of the behavior of
$\rho_E$ near 0, it will be profitable to find a simple majorant
sufficiently close to $\rho_E$.

\vskip 3mm

\noindent \textbf{Definition 2.3.} A function $\rho\! :\!( 0,1]
\rightarrow (0,+\infty )$  {\it satisfies conditions of a fit
majorant} if

\vskip 1mm 1) $t\longmapsto\rho (e^t)$ \ is an increasing,  $
\mathcal{C}^1((-\infty,0])$ strictly convex function,

2) $\rho '(1)\geq 1$ \ and \ $\lim\limits_{r\rightarrow 0+}\rho
(r)=0$.
 \vskip 1mm

Note that for every function $\rho$ satisfying the above
conditions we have $\lim_{t\rightarrow -\infty}\psi (t)=0$ and
$\psi^{-1}((0,\rho'(1)])=(-\infty ,0]$ where $\psi (t):=(\rho
(e^t))'$.

 \vskip 3mm

\noindent \textbf{Example 2.4.} The following functions satisfy
conditions of a fit majorant:

\vskip 1mm i) \ $\rho (r)=Ar^\gamma$ \ for \ $\gamma\in (0,1]$ and
$A\geq 1/\gamma$,

ii) \ $\rho (r)=\frac{1}{s}(1/\log (e/r))^s$ \ for  $s>0$,

iii) \ $\rho (r)=Ar^\sigma \left(\log
(1/r)+\frac{2}{\sigma}\right)$ \ for  $\sigma \in (0,1]$ and $
A\geq 1$.

\vskip 3mm

We are interested in seeing how the H\"older continuity of the
pluricomplex Green's function $V_E$ is connected with Markov-type
inequalities for polynomials on $E$. The question will be made
more precise by the next definition:

\vskip 3mm

\noindent \textbf{Definition 2.5.} Let $\gamma\in (0,1], \ B>0$. A
compact set $E\subset \mathbb{C }^N$ admits the {\textit{H\"older
continuity property of the pluricomplex Green's function $V_E$}}
($E\in HCP(\gamma, B$) in short) if for every $z\in \mathbb{C}^N$
\be V_E(z) \ \le \ B \ {\rm dist} (z,E)^{\: \gamma}. \label{HCP}
\ee

\vskip 3mm

It seems appropriate to mention here five equivalents for this
property.

vskip 3mm

\noindent \textbf{Proposition 2.6.} \textit{If $E$ is a compact
subset of $\mathbb{C}^N$ and $\gamma\in (0,1]$ then the following
statements are equivalent:}

\vskip 1mm

$(i) \ \ \ \, \exists \: B_1\ge 1 \ \ \ E\in HCP(\gamma, B_1),$

\vskip 1mm

$(ii) \ \, \ \exists \: B_2\ge 1 \ \ \ \rho_E(r) \le B_2\,
r^\gamma \ \ \ \ \ {\it for} \ \ r\ge 0,$

\vskip 1mm

$(iii) \ \, \exists \: B_3\ge 1 \ \ \ |\rho_E(r)-\rho_E(s)|\le
B_3\, |r-s|^\gamma \ \ \ \ \ {\it for} \ \ r,s\ge 0,$

\vskip 1mm

$(iv) \ \, \, \exists \: B_4\ge 1 \ \ \ \Phi_E(z) \le 1+ B_4 \
{\rm dist}\, (z,E)^{\: \gamma} \ \ \ \ {\it for} \ \ \ z\in
\mathbb{C}^N, \ \ {\rm dist}\, (z,E)\le 1,$

\vskip 1mm

$(v) \ \ \ \exists \: B_5\ge 1 \ \ \ | V_E(z) - V_E(w) | \ \le \
B_5 \ \|z-w\|_2^{\: \gamma} \ \ \ \ {\it for} \ \ \ z,w\in
\mathbb{C}^N,$

\vskip 1mm

$(vi) \  \ \forall \: R>0 \ \ \exists \: B_6\ge 1 \ \ \ |
\Phi_E(z) - \Phi_E(w) | \ \le \ B_6 \ \|z-w\|_2^{\: \gamma} \ \ \
\ {\it for} \ \ \ z,w\in E_R:=\{z\in \mathbb{C}^N\: : \: {\rm
dist}\, (z,E)\le R\}$.

\vskip 1mm

\noindent {\it Moreover, in the equivalences $(i) \Leftrightarrow
(ii) \Leftrightarrow (iii) \Leftrightarrow (v)$ we have
$B_1=B_2=B_3=B_5$.}

\vskip 2mm

If $E\in HCP(\gamma ,B)$ then $E$ is $L$-regular and therefore
$C(E)>0$. However, a lower bound for $C(E)$ in terms of the
constants $\gamma ,B$ was not known. In this paper we simply solve
this problem (see Th.2.12).

A close inspection of the proof of \cite[Th.3.5]{[15]} and use of
Stirling's approximation lead us to

\vskip 3mm

\noindent \textbf{Proposition 2.7.} {\it If $E\subset \mathbb{C}$
and there exists $M_k=M_k(E)$ independent of $n$ and a polynomial
$P$ of degree at most $n$ such that \be \| P^{(k)} \|_E \ \leq \
M_k \ n^{mk} \: \|P\|_E \label{M_k} \ee then \ $ M_k \ \ge \
B^{k}/[({k}!)^{m-1}]$ \ for certain absolute constant $B>0$.}

\vskip 3mm

This fact was the motivation for concentrating on the following
generalization of the V.Markov inequality.

Let $\mathbb{N}=\{1,2,...\}$ and
$\mathbb{N}_0=\mathbb{N}\cup\{0\}$.

\vskip 2mm

 \noindent
\textbf{Definition 2.8.} Fix $m\ge 1, \ M>0$.  A compact set
$E\subset \mathbb{C }^N$ admits the {\it V.Markov inequality }
($E\in VMI(m, M$) in short) if for every $\alpha\in
\mathbb{N}_0^N$, $P\in {\cal{P}}(\mathbb{C}^N)$  \be
 \|D^\alpha P \|_E \
\leq \ M^{|\alpha|} \ \frac{({\rm deg}\,P)^{m{|\alpha|}}}{
({|\alpha|}!)^{m-1}} \: \|P\|_E \label{VMI} \ee where
$|\alpha|=\alpha_1+\ldots+\alpha_N$, $D^\alpha =
\frac{\partial^{|\alpha|}}{\partial z_1^{\alpha_1}\ldots\,
\partial z_N^{\alpha_N}}$ \ for $\alpha=(\alpha_1,\ldots,\alpha_N)$.

\vskip 2mm

In other words, (\ref{VMI}) is a version of inequality (\ref{M_k})
(and also its analogue in higher dimensional space) with the
strongest possible constants $M_k$ (compare with \cite{bbcm}).

\vskip 3mm

\noindent \textbf{Example 2.9.} The simplest example of a set
admitting V.Markov inequality is the unit disc
$\overline{\mathbb{D}}$ in the complex plane. By Bernstein
inequality, for every polynomial $P$ of degree not greater than n
we have \ $\|P^{(k)} \|_{\overline{\mathbb{D}}} \leq
\frac{n!}{(n-k)!} \ \| P\|_{\overline{\mathbb{D}}}$ \ and thus \
$\overline{\mathbb{D}}\in VMI(1,1)$. In the multidimensional space
for a polydisc $D(a,r) = \{z\in \mathbb{C}^N \, : \, |z_1-a_1|\le
r_1,\ldots, |z_N-a_N|\le r_N\}$ of polyradius $r=(r_1,\ldots,
r_N)\in (0,+\infty)^N$ we get (see \cite[Example 2.2]{[6]}) \[
\|D^\alpha P \|_{D(a,r)} \ \leq \ \frac{\nu !}{(\nu - \alpha)!\:
r^\alpha} \ \|P\|_{D(a,r)} \] whenever $P$ is a polynomial of $N$
variables $z_1,\ldots,z_N$ of degree at most $\nu_1$ in $z_1$, ...
, $\nu_N$ in $z_N$. As usual, $\nu !=\nu_1!\ldots \nu_N!$ and
$r^\alpha=r_1^{\alpha_1}\!\ldots r_N^{\alpha_N}$. Hence \be
\|D^\alpha P \|_{D(a,r)} \ \le \ {r^{-\alpha}} \
\nu_1^{\alpha_1}\ldots\nu_N^{\alpha_N} \ \|P\|_{D(a,r)} \ \le \
r^{-\alpha} \ ({\rm deg}\,P)^{|\alpha|} \ \|P\|_{D(a,r)} \label{3}
\ee for any polynomial $P\in {\cal{P}}(\mathbb{C}^N)$. Therefore
${D(a,r)}\in VMI(1, \max\limits_{j} 1/r_j)$.

\vskip 2mm

\noindent \textbf{Example 2.10.} Due to the classical inequality
proved by V.Markov (see (\ref{VM})), we have $[-1,1] \!\in\!
VMI(2,1)$. If $E\!=\![a_1,b_1]\! \times\! \ldots\! \times\!
[a_N,b_N]\! \subset\! \mathbb{R}^N\! \subset\! \mathbb{R}^N \!+\!
i\mathbb{R}^N\! = \mathbb{C}^N$ then for every polynomial $P$ of
degree at most $\nu_1$ in $z_1$, \dots , $\nu_N$ in $z_N$  we have
(see \cite[Example~2.2]{[6]})
\[  \|D^\alpha P
\|_E \ \leq \ \frac{2^{|\alpha|}}{(b-a)^\alpha} \
T_{\nu_1}^{(\alpha_1)}(1) \cdot \ldots \cdot
T_{\nu_N}^{(\alpha_N)}(1) \ \|P\|_E \ \ \le \
\frac{2^{{|\alpha|}}}{(b-a)^\alpha} \cdot
\frac{\nu^{2\alpha}}{\alpha !} \ \|P\|_E \] with $a=(a_1,\ldots
,a_N), \ b=(b_1,\ldots ,b_N), \ \nu=(\nu_1,\ldots,\nu_N)$.   Since
$N^{|\alpha|}\alpha!\ge |\alpha|!$, we obtain $E\in
VMI\left(2,2N\max\limits_{j} 1/(b_j-a_j)\right)$.

\vskip 3mm

It is evident that $VMI(m,M) \Rightarrow AMI(m, M\sqrt{N})$. On
the other hand, property (\ref{AMI}) easily implies that
\[ \|D^\alpha P \|_E \ \leq \ M^{|\alpha|} \ \left( \frac{n!}{(n-|\alpha|)!}\right)^{m} \: \|P\|_E \ \leq \ M^{|\alpha|} \ n^{m{|\alpha|}} \: \|P\|_E \] for
any $\alpha\in \mathbb{N}_0^N$, $P\in {\cal{P}}(\mathbb{C}^N)$ of
degree at most $n$.

\vskip 3mm

\noindent \textbf{Remark 2.11.}  If $E\in AMI(m_1, M_1)$ and if we
fix an arbitrary $\delta\in (0,1)$ then for every polynomial $P$
of degree at most $n$ and for all $|\alpha| \le n^\delta$,
inequality (\ref{VMI}) holds with $m=\frac{m_1-\delta}{1-\delta}$
and $M=M_1$. Indeed,
\[ \|D^\alpha P\|_E \ \le \ (M_1\, n^{m_1})^{|\alpha|} \|P\|_E\]
and  \[ n^{m_1|\alpha|} = \left(
\frac{n^{m_1}|\alpha|^{m-1}}{|\alpha|^{m-1}} \right)^{|\alpha|}
\le
\frac{(n^{m_1+(m-1)\delta})^{|\alpha|}}{|\alpha|^{(m-1)|\alpha|}}
\le \frac{n^{m|\alpha|}}{|\alpha|!^{m-1}}.\] By the above, in the
particular case of $m_1=1$, we get $AMI(1, M_1) \Leftrightarrow
VMI(1,M)$.

\vskip 3mm

In the general case, we do not know whether or not the V.Markov
inequality is equivalent to that of A.Markov. However, we can show
that the H\"older continuity property is equivalent to
(\ref{VMI}).

\vskip 3mm

\noindent \textbf{Theorem 2.12 (Main theorem).} \textit{If $E$ is
a compact subset of $\mathbb{C}^N$, $0< \gamma\le 1 \le m$,
$B,M>0$ then \be E\in HCP(\gamma, B) \ \ \Longrightarrow \ \ E\in
VMI\left(m,M \right) \ \ \ {with} \ \ m=1/\gamma,\ \ M=\sqrt{N}\:
(B\gamma e)^{1/\gamma} \label{hcpvmi}\ee \be E\in VMI(m, M) \ \
\Longrightarrow \ \ E\in HCP\left(\gamma, B\right) \ \ \ {with} \
\ \gamma=1/m, \ \ B=M^{\gamma}N^{\gamma} m . \ \ \ \ \ \ \ \ \
\label{vmihcp}\ee}\noindent{\it Moreover, if $ E\in VMI\left(m,M
\right)$, then $C(E)\geq e^{-m} \frac{1}{NM}$. Hence, if $E\in
HCP\left(\gamma, B\right)$, then $C(E)\geq \left(N^{3/2}(B\gamma
e^2)^{1/\gamma}\right)^{-1}$.}

 \vskip 1mm

As a consequence of the above theorem, the well known open problem
concerning the conjectured implication $AMI\Rightarrow HCP$ is
equivalent to a new question of whether $AMI$ implies $VMI$. The
first problem regards the properties related to the notions in two
different fields: the pluricomplex Green's function and
polynomials, whereas the new question is formulated only in terms
of derivatives of polynomials.

\vskip 2mm

Due to the above theorem, we can give new, somewhat unexpected
equivalents to the H\"older continuity property of the
pluricomplex Green's function:

\vskip 3mm

\noindent \textbf{Corollary 2.13.} \textit{If $E$ is a compact
subset of $\mathbb{C}^N$ and $\gamma\in (0,1]$ then the following
conditions are equivalent:}

\vskip 2mm

$(i) \ \ \ E\in HCP(\gamma, B_1) \ \ \ {\it with \ some} \ \ \
B_1\ge 1,$

\vskip 2mm

$(ii) \ \, \ \exists \: B_2>0 \ \ \forall \: z_0\in E \ \ \forall
j\in\{1,\ldots,N\} \ \ \forall \: \zeta\in \mathbb{C}$ \ {\it such
that} \ $|\zeta|\le 1$ {\it we have} \[ V_E(z_0+\zeta e_j) \ \le \
B_2 \: |\zeta| ^\gamma,\]

\vskip 2mm

$(iii) \ \, \exists \: M_3 > 0 \ \ \forall \: j \in \{1,\ldots,N\}
\ \ \forall \: P \in \mathcal{P} (\mathbb{C}^N) \ \ \forall \:
k\in\mathbb{N} $ {\it we have} \[ \|D^{ke_j} P \|_E \ \leq \
M_3^{k} \ \frac{({\rm deg}\,P)^{k/\gamma}}{ k\,!^{\,
\frac1\gamma-1}} \ \|P\|_E,\]

\vskip 2mm

\noindent {\it where $e_1,\ldots,e_N$ are the canonical vectors in
$\mathbb{C}^N$: $e_j=(0,\ldots,0,1,0,\ldots,0)$ with the value 1
in the jth entry.}

\vskip 3mm

It seems to be rather surprising that condition $(ii)$ in Cor.2.13
that holds only in $N$ canonical directions, is sufficient to
guarantee the H\"older continuity property of $V_E$ in all
directions.

We can generalize the main theorem  to the case where $\rho_E(r)$
and $k$th derivatives of polynomials have bounds related to fit
majorants with some additional properties (Th.4.2). We show that
the required properties are satisfied for functions given in
Example 2.4 i and ii (Th.4.4, 4.5).

For the compact subsets of $\mathbb{R}^N \subset \mathbb{R}^N +i
\mathbb{R}^N = \mathbb{C}^N$ we prove that if inequality
(\ref{HCP}) holds for $x\in \mathbb{R}^N$ then so is for all $z\in
\mathbb{C}^N$ (Cor.5.4). As a consequence, we obtain that $E\in
HCP\left(1/2, {\mathcal{B}}/{\sqrt{C(E)}}\right)$ for any convex
body in $\mathbb{R}^N$, where $C(E)$ is the L-capacity of $E$ and
$\mathcal{B}$ is an absolute constant independent of $E$ and even
of $N$ (Example 5.7). Moreover, we prove that every set $E\subset
\mathbb{R}^N$ uniformly polynomially cuspidal in direction $v$
with exponent $s$, has the following property: $V_E(x+\zeta v)\leq
B|\zeta|^{1/(2s)}$ for $x\in E,\ |\zeta |\leq r_0$ (Th.5.10).
Hence we deduce from Cor.2.13 that every UPC compact subset of $
\mathbb{R}^N$ admits HCP and thus V.Markov inequality (Cor.5.11).
In this way, we obtain a wide class of sets that have such a
property. This is the first essential generalization of V.Markov's
result from the end of XIX century.

As another application of the main theorem, we can prove HCP for
disconnected sets. Prop.6.1 regards some onion type sets in the
complex plane that may not satisfy the assumptions of the
Carleson-Totik criterion. These sets are particularly interesting
in view of certain properties of compacts admitting so-called
local Markov's inequality (see L.Bialas-Ciez and R.Eggink, {\it
Equivalence of the global and local Markov inequalities in the
complex plane}, in preparation). The second example of such an
application of Th.2.12 concerns some compact sets consisting of
infinitely many pairwise disjoint subsets of $\mathbb{C}^N$
(Prop.6.2).

\vskip 8mm

{\bf  3. Proofs of the main results}

\vskip 3mm

To prove Prop.2.6 we need the following

\vskip 3mm

\noindent \textbf{Lemma 3.1.} \textit{If $E\subset\mathbb{C}^N$ is
a compact L-regular set then} \be a) \ \ \ |V_E(z)-V_E(w)| \ \le \
\rho_E(\|z-w\|_2) \ \ \ \ \ \ \ \ \ \ {for} \ \ \ z,w\in
\mathbb{C}^N, \ \ \ \ \ \ \
 \ \ \ \ \ \ \ \ \ \ \ \ \ \ \
 \
 \ \label{rho1} \ee
\[ b) \ \ \ |\rho_E(r)-\rho_E(s)| \ \le \ \rho_E(|r-s|) \ \ \
{for} \ \ r,s\ge 0. \ \ \ \ \ \ \
 \ \ \ \ \ \ \ \
 \ \ \ \ \ \ \ \ \ \ \ \ \ \ \
 \ \ \ \ \ \ \ \
  \ \ \ \
 \
 \ \]

\vskip 1mm

\noindent \textit{Proof of Lemma 3.1.} Modifying an argument due
to B\l ocki (see \cite[Prop.3.5]{[21]}), consider the function
$u_\zeta(z)\!:=\!V_E(z\!+\!\zeta)\!-\! \rho_E(\|\zeta\|_2)$ for
$z\in \mathbb{C}^N$. We have $u_\zeta \in {\cal{L}}_N$ for any
fixed $\zeta\in \mathbb{C}^N$. Moreover, $$V_E(z+\zeta) \ \le \
\|V_E\|_{E_{{\rm dist}\, (z+\zeta,E)}} \ \le \
\|V_E\|_{E_{\|\zeta\|_2}} \ = \  \rho_E(\|\zeta\|_2)$$ whenever
$z\in E$. Hence $u_\zeta \le 0$ on $E$ and, by the definition of
the pluricomplex Green's function, $V_E(z)\ge u_\zeta(z)$ for all
$z \in \mathbb{C}^N$ and we can easily obtain statement $a)$.

To prove b), fix $r\ge s\ge 0$. We can take $w\in E_r$ such that
$V_E(w)=\rho_E(r)$. Choose $z\in E_s$ at distance $r-s$ from the
point $w$. By inequality $(\ref{rho1})$, we have
\[ 0\le \rho_E(r)-\rho_E(s) \ \le \ V_E(w) -
V_E(z) \ \le \ \rho _E (\|w-z\|_2)=\rho_E(r-s)\] and we get
property b). \hfill $\square$

\vskip 5mm

\noindent \textit{Proof of Proposition 2.6.} The implication
$(i)\Rightarrow (ii)$ is an easy consequence of the definition of
$\rho_E$.

Lemma 3.1 immediately implies $(iii)$ whenever we assume $(ii)$
and we take into account Remark 2.2.

If we consider $s=0$ and $r={\rm dist}\,(z,E)$ for a fixed $z\in
\mathbb{C}^N$,  we obtain $(i)$ from $(iii)$.

The equivalence $(i) \Leftrightarrow (iv)$ is an easy consequence
of (\ref{V_E}) and the elementary inequalities: $1+x\le e^x$ for
$x\ge 0$ and $e^{bt} \le 1+(e^b-1)t$ \ for $t\in [0,1], \ b>0$.

By B\l ocki's argument mentioned in the proof of Lemma 3.1,
property $(i)$ implies the H\"older continuity of the pluricomplex
Green's function $V_E$ in the whole space, i.e. condition $(v)$.

To prove $(v)\Rightarrow(vi)$, it is sufficient to apply formulas
(\ref{V_E}) and the fact that $V_E$ is lower semicontinuous.
Indeed, for $z,w\in E_R$
\[ | \Phi_E(z) - \Phi_E(w) | = |\exp (V_E(z))- \exp (V_E(w))| \] \[\le
| V_E(z)- V_E(w)| \: \exp ( \max \{ V_E(z), V_E(w) \} )  \le B_5 \
\|z-w\|_2^{\: \gamma} \: \exp ( \|V_E\|_{E_R}).\]

The evident implication $(vi) \Rightarrow (iv)$ finishes the
proof. \hfill $\square$

\vskip 5mm

\noindent \textit{Proof of Theorem 2.12.} To show the first
implication, consider an arbitrary polynomial $P\in
{\cal{P}}(\mathbb{C}^N)$ of degree at most $n$ and $\alpha\in
\mathbb{N}_0^N$. By Cauchy's integral formula and the
Bernstein-Walsh-Siciak inequality, for fixed $z\!=\!(z_1,\ldots,
z_N)\!\in \!E$, $r\!\in \!(0,1]$ we can obtain
\[ |D^\alpha P(z)| \le \frac{\alpha !}{(r/\sqrt{N})^{|\alpha|}} \:
\| P\|_{D(z,r/\sqrt{N})} \le \frac{\sqrt{N}^{|\alpha|} \alpha
!}{r^{|\alpha|}} \: \| P\|_E \ \exp\left(n
\|V_E\|_{D(z,r/\sqrt{N})}\right)\] where $D(z,r/\sqrt{N})=\{w\in
\mathbb{C}^N \, : \, |w_1-z_1|\le r/\sqrt{N},\ldots, |w_N-z_N|\le
r/\sqrt{N}\}$. From (\ref{HCP}) we have \[ |D^\alpha P(z)|  \le
\frac{\sqrt{N}^{|\alpha|} \alpha !}{r^{|\alpha|}} \: \| P\|_E \
\exp\left(n \:B{r^\gamma}\right)\] and for
$r=\left({|\alpha|}/(B\gamma n)\right)^m$,  $m=1/\gamma$ we get
\[ |D^\alpha P(z)| \le \frac{\sqrt{N}^{|\alpha|} \alpha
!}{|\alpha|^{|\alpha|m}} \: n^{m|\alpha|}\: \| P\|_E \ (B\gamma
e)^{ |\alpha|/\gamma } \le \left(\sqrt{N}(B\gamma
e)^{1/\gamma}\right)^{|\alpha|} \frac{|\alpha| !}{(|\alpha|!)^{m}}
\: n^{m|\alpha|}\: \| P\|_E ,
\] and (\ref{hcpvmi}) is proved.

We now proceed to show implication (\ref{vmihcp}). For this
purpose, observe that from (\ref{V_E}), it is sufficient to prove
\be |P(z)| \ \le \ \|P\|_E \ \exp\left(M^{\gamma}N^{\gamma} m
nr^\gamma\right) \label{wyst} \ee for any polynomial $P\in
{\cal{P}}(\mathbb{C}^N)$ of degree at most $n$ and $z\in
\mathbb{C}^N\setminus E$ such that dist$(z,E)=r$. By Taylor's
formula, we have
\[ |P(z)| \ \leq \ \sum\limits_{|\alpha|\le n} \frac{1}{\alpha !} \
|D^\alpha P(w)| \ r^{|\alpha|} \ \le \ \sum_{k=0}^n
\sum\limits_{|\alpha|=k} \frac{1}{\alpha !} \ r^{|\alpha|} \:
\|D^\alpha P\|_E \] whenever $w\in E$ and dist$(z,E)=\|z-w\|_2$.\
 From (\ref{VMI}) the above inequality gives
\[ |P(z)| \ \leq \ \|P\|_E \ \sum_{k=0}^n M^{k} \, \frac{n^{m{k}}}{
({k}!)^{m-1}} \ r^{k} \sum\limits_{|\alpha|=k} \frac{1}{\alpha !}
. \] Since $\sum_{|\alpha|=k} {1}/{\alpha !} \ = \ {N^k}/{k!}$, \
for $\gamma=1/m$ we have
\[ |P(z)| \ \leq \ \|P\|_E \ \sum_{k=0}^n M^{k}N^k \ \frac{n^{m{k}}}{
({k}!)^{m}} \ r^{k}  \ \leq \ \|P\|_E \ \sum_{k=0}^\infty \left[
\frac{\left( M^{\gamma}N^{\gamma}r^{\gamma} n\right)^k}{
{k}!}\right]^{m} \] \[ = \ \|P\|_E \ {\cal{G}}_m \left(
M^{\gamma}N^{\gamma}r^{\gamma} n\right)  \] where
\[ {\cal{G}}_m(x):=\sum\limits_{k=0}^\infty \left( \frac{x^k}{k!}
\right)^{\!m} = \left\| \left( \frac{x^k}{ {k}!}\right)_{\! k\in
\mathbb{N}_0} \right\|_m^{\, m} \] and $\|\cdot \|_m$ is the usual
norm in the space ${{l}}_m$. As $\|\cdot \|_m \le \|\cdot \|_1$,
we have \[{\cal{G}}_m(x) = \left\| \left( \frac{x^k}{
{k}!}\right)_{\! k\in \mathbb{N}_0} \right\|_m^{\, m} \le \left\|
\left( \frac{x^k}{ {k}!}\right)_{\! k\in \mathbb{N}_0}
\right\|_1^{\, m} = e^{xm} \] and \[ |P(z)| \ \leq \  \|P\|_E \
\exp \left( M^{\gamma}N^{\gamma}r^{\gamma} m n\right),  \] which
gives (\ref{wyst}).

If $w\in E$ and $\zeta\in\mathbb{C}^N,\ \|\zeta\|_2=r$, then in a
similar way as above we get
$$|P(w+\zeta )| \ \leq \ \|P\|_E \ \sum_{k=0}^n N^{k}M ^k \ \frac{n^{m{k}}}{
({k}!)^{m}} \ r^{k}.$$ In the case of $NMr\geq 1$ we obtain
$$|P(w+\zeta )| \ \leq \ \|P\|_E
(NMr)^n\left(\sum_{k=0}^n\frac{n^k}{k!}\right)^m\leq \|P\|_E
(NMr)^ne^{nm}.$$ Thus $\rho_E(r)\leq \log (NMe^m)+\log r$ \ for
$r\geq e^{-m}\frac{1}{NM}$ \ and consequently (see Remark 2.2),
$$-\log C(E)=\lim\limits_{r\rightarrow\infty}(\rho_E(r)-\log
r)\leq \log (NMe^m)$$
 and the proof is completed. \hfill $\square$

\vskip 3mm

\noindent \textit{Proof of Corollary 2.13.} First, we prove
$(ii)\Rightarrow(iii)$. Put $p_j(\zeta)=P(z_0+\zeta e_j)$ for
$\zeta\in \mathbb{C}$, $z_0\in E$ and for a fixed polynomial $P\in
{\cal{P}}(\mathbb{C}^N)$ of degree at most $n$. Obviously, $
|D^{ke_j} P(z_0)| = |p_j^{(k)}(0)|$ and by Cauchy's integral
formula,
\[ |D^{ke_j}P(z_0)| \le \frac{k!}{r^k} \max\{|p_j(\zeta)| \: : \:
|\zeta|=r\} \le \frac{k!}{r^k} \|P\|_E\max\{\exp(n V_E (z_0+\zeta
e_j)) \: : \: |\zeta|=r\} \] the last inequality being  a
consequence of (\ref{V_E}). If $r=\left({k}/{n}\right)^{1/\gamma}$
then from $(ii)$ we get
\[ \|D^{ke_j}P\|_E \le k! \, \left(\frac{n}{k}\right)^\frac{k}\gamma  e^{B_2\, k} \: \|P\|_E \le M_3^{k} \ \frac{n^{k/\gamma}}{ k\,!^{\,
\frac1\gamma-1}} \ \|P\|_E\] with some positive constant $M_3$,
and $(iii)$ is proved.

In view of Th.2.12, to show $(iii)\Rightarrow(i)$ it is sufficient
to prove that $(iii)$ implies inequality (\ref{VMI}). Fix
$\mathbb{N}_0^N\ni \alpha = \alpha_1 e_1 + \ldots +\alpha_N e_N $.
If $P$ is a polynomial of degree $n_j$ in $z_j$ where
$z=(z_1,\ldots,z_N)$, we have \[ \| D^\alpha P\|_E \ \le \
M_3^{\alpha_1} \, \frac{n_1^{\alpha_1/\gamma}}{ \alpha_1!^{\,
\frac1\gamma-1}} \ \| D^{\alpha - \alpha_1 e_1} P\|_E \ \le \
\ldots\]
\[\le \ M_3^{|\alpha|} \, \frac{n_1^{\alpha_1/\gamma}\ldots
n_N^{\alpha_N/\gamma}}{ \alpha\,!^{\, \frac1\gamma-1}} \: \| P\|_E
\ \le \ M_3^{|\alpha|}\, N^{|\alpha|} \ \frac{({\rm
deg}\,P)^{|\alpha|/\gamma}}{ |\alpha|\,!^{\, \frac1\gamma-1}} \
\|P\|_E\] since $N^{|\alpha|}\alpha!\ge |\alpha|!$, and $(i)$
follows.

 The
implication $(i)\Rightarrow(ii)$ is obvious and the proof is
completed. \hfill $\square$ \vskip 3mm

\noindent {\bf Remark 3.2.} It follows from the proof of Cor.2.13
that we can replace condition $V_E(x\!+\!\zeta e_j)\!\leq\!
C_2|\zeta |^\gamma$ for $|\zeta|\leq 1$ by the same but only for
$|\zeta |\leq r_0\leq 1$. Indeed, in the proof of implication
$(ii)\Rightarrow (iii)$ it is sufficient to put
$r=r_0(k/n)^{1/\gamma}$. We shall use this remark later.

\bigskip

{\bf 4. Majorants of $V_E$ and a bound of $k$th derivative of
polynomials} \vskip 3mm

In this section we shall present a generalization of properties
related to HCP. To do this, we need the following basic
definition.

\vskip 3mm
 \noindent \textbf{Definition 4.1.} Let $\rho$ satisfy conditions of a fit
majorant.

\noindent a) We say that $\rho$ is {\it m-bounded} if  for all
$c\in (0,1]$ and $ r\in [0,1]$ \be
\limsup\limits_{n\rightarrow\infty} \log
\left(1+\sum\limits_{k=1}^n\exp\left(-k\psi^{-1}\left( ck/n\right)
+n\rho \left(\exp\psi^{-1}\left(
ck/n\right)\right)\right)r^k\right)^{1/n}\leq \frac{A}{c}\rho (r),
\label{m}\ee where $\psi (t)=(\rho (e^t))'$ and $A$ is a constant
independent of $c$ and $r$.

\noindent b) We say that $\rho$ is {\it mb-bounded} if for all
$s\in\mathbb{N}$ there exist constants $C_s\geq 0,c_s\in (0,1]$
such that for all $1\leq k_1+\dots +k_s\leq n$
$$-\sum\limits_{j=1}^{s} \ \!\!\! 'k_j\psi^{-1}\left(
\frac{k_j}{n}\right)+n\sum\limits_{j=1}^{s} \ \!\!\! '\rho\left(
\exp\psi^{-1}\left( \frac{k_j}{n}\right)\right)$$ $$ \leq
-\sum\limits_{j=1}^sk_j\psi^{-1}\left(c_s
\frac{1}{n}\sum\limits_{j=1}^sk_j\right)
+n\rho\left(\exp\psi^{-1}\left(
c_s\frac{1}{n}\sum\limits_{j=1}^sk_j\right)\right)+
C_s\sum\limits_{j=1}^sk_j,$$ where $\sum '$ means that we consider
only $k_j$ such that $k_j\neq 0$.

\noindent c) We say that $\rho$ is {\it doubly bounded} if \
$\limsup\limits_{r\rightarrow 0+}\rho (2r)/\rho (r)< +\infty$.

\vskip 1mm

 \noindent \textbf{Theorem 4.2.} {\it Let $E$ be a compact
subset of $ \mathbb{C}^N$ and let $\rho$ satisfy conditions of a
fit majorant. Put $\psi (t)=(\rho (e^t))'$.

\noindent a) Assume that \ $ \rho_E(r)\leq \rho (r)$ \ for $ r\in
(0,1]$. For all $n\geq 1,\ \alpha\in \mathbb{N}_0^N,\ 1\leq
|\alpha |\leq n$ and for an arbitrary $P\in \mathcal{P}(
\mathbb{C}^N)$ of degree not greater than $n$ we have
$$||D^\alpha P||_E\leq \alpha ! N^{|\alpha |/2}\exp\left(-|\alpha|\psi^{-1}\left(
|\alpha|/n\right) +n\rho \left(\exp\psi^{-1}\left(
|\alpha|/n\right)\right)\right) ||P||_E.$$ \vskip 2mm

\noindent b)  If $\rho$ is {\it m}-bounded (with a constant $A$)
and there exist constants $C>0,c\in (0,1]$ such that for all
polynomials $P$ of degree not greater than $n$ and $1\leq
|\alpha|\leq n$
$$||D^\alpha P||_E\leq \alpha! C^{|\alpha|}\exp\left(-|\alpha |\psi^{-1}\left(
c|\alpha |/n\right) +n\rho \left(\exp\psi^{-1}\left( c |\alpha
|/n\right)\right)\right) ||P||_E$$ then we have the inequalities
$$\rho_E(r)\leq \frac{A}{c}\rho ( Cr),\ \ \ \ for \ \ 0\leq r\leq
1/C,$$ $$\rho_E(r)\leq \frac{A}{c}\rho (1)+\log (Cr),\ \ \ \ for \
\  r\geq 1/C$$ and $C(E)\geq \exp (-A\rho (1)/c)\, /C$. \
Moreover, if $\rho$ is doubly bounded, then there exists a
constant $A'$ such that $\rho_E(r)\leq A'\rho (r),\ r\in [0,1]$.
\vskip 2mm

 \noindent c) If $\rho$ is {\it m}-, {\it mb}- and doubly bounded
and $V_E(x+\zeta e_j)\leq \rho (|\zeta|),\ x\in E,\zeta\in
\overline{\mathbb{D}}$, then there exists a constant $B$ such that
$V_E(x+\zeta )\leq B\rho (||\zeta ||_2)$ for $ x\in E,
\zeta\in\mathbb{C}^N,\|\zeta\|_2\leq 1$.}

\vskip 3mm

\noindent {\it Proof.} a) By Cauchy's inequality, we have
$$|D^\alpha P(x)|/\alpha!\leq  \inf\limits_{r>0}\{ r^{-|\alpha
|}\sup\limits_{|\zeta_j|\leq r}\exp (nV_E(x+(\zeta_1,\dots
,\zeta_N)))\}$$
$$\leq  \inf\limits_{r>0}\{ r^{-|\alpha
|}\sup\limits_{||\zeta||_2\leq \sqrt{N} r}\exp (nV_E(x+\zeta
))\}\leq \inf\limits_{r>0} r^{-|\alpha |}\exp (n\rho_E(\sqrt{N}r)
$$ $$=N^{|\alpha|/2}\inf\limits_{r>0} r^{-|\alpha |}\exp (n\rho_E(r).
$$ Hence
$$ ||D^\alpha P(x)||_E/\alpha! \leq N^{|\alpha |/2}\inf\limits_{r\in (0,1]}
r^{-|\alpha|}e^{n\rho (r)}=N^{|\alpha |/2} \inf\limits_{t\leq
0}\exp \left(-|\alpha |t+n\rho (e^t)\right)$$
$$=N^{|\alpha |/2}\exp\left(\inf\limits_{t\leq 0} \left(-|\alpha |t+n\rho
(e^t)\right)\right).$$ Since $f(t)=-|\alpha |t+n\rho (e^t)$ \ is a
strictly convex function and $\rho'(1)\geq 1$ we get
$$\inf\limits_{t\leq 0} \left(-|\alpha |t+n\rho
(e^t)\right)=-|\alpha|\psi^{-1}(|\alpha |/n)+n\rho(\exp
\psi^{-1}(|\alpha |/n)),$$ and assertion a) follows. \vskip 1mm

\noindent b)  By  Taylor's theorem applied to a polynomial $P$ of
degree $n\geq 1$, $||P||_E=1$, for $x\in
E,\zeta\in\mathbb{C}^N,||\zeta||_2=r$ we can write

$$|P(x+\zeta )|\leq |P(x)|+\sum\limits_{0<|\alpha|\leq n} \frac{|D^\alpha
P(x)|}{\alpha!}\: ||\zeta ||_2^{|\alpha|}$$ $$\leq
1+\sum\limits_{k=1}^n\sum\limits_{|\alpha
|=k}C^{|\alpha|}\exp\left(-|\alpha|\psi^{-1}\left(
c\frac{|\alpha|}{n}\right) +n\rho \left(\exp\psi^{-1}\left(
c\frac{|\alpha|}{n}\right)\right)\right) r^{|\alpha|}$$
$$=1+\sum\limits_{k=1}^n\binom{N+k-1}{k}C^{k}\exp\left(-k\psi^{-1}\left(
ck/n\right) +n\rho \left(\exp\psi^{-1}\left(
ck/n\right)\right)\right)r^{k}$$
$$\leq 1+\sum\limits_{k=1}^n\frac{(N+k-1)^{N-1}}{(N-1)!}\exp\left(-k\psi^{-1}\left(
ck/n\right) +n\rho \left(\exp\psi^{-1}\left(
ck/n\right)\right)\right) (Cr)^{k}$$
$$\leq \frac{(N+n-1)^{N-1}}{(N-1)!}\left(1+\sum\limits_{k=1}^ne^{(N-1)k}\exp\left(-k\psi^{-1}\left(
ck/n\right) +n\rho \left(\exp\psi^{-1}\left(
ck/n\right)\right)\right) (Cr)^{k}\right).$$
 Replacing $P$
by $P^m$ we get $$ \log |P(x+\zeta )|^{1/n}\leq \log\left(
(mn+N-1)^{N-1}/(N-1)!\right)^{1/mn} $$
$$ +\log \left(1+\sum\limits_{k=1}^{mn}\exp\left(-k\psi^{-1}\left(
ck/mn\right) +mn\rho \left(\exp\psi^{-1}\left(
ck/mn\right)\right)\right) (Cr)^{k}\right)^{1/mn}.$$
 Hence, from (\ref{m})  we
obtain for $Cr\leq 1$ the inequality
$$\log |P(x+\zeta )|^{1/n}\leq \frac{A}{c}\rho
(Cr)$$ and consequently,
$$\rho_E(r)\leq \frac{A}{c}\rho (Cr),\ r\leq 1/C.$$ If $Cr\!\!\geq\!\! 1$
then $(Cr)^k\!\!\leq\!\! (Cr)^n$ for $k\!\leq \!n$ and a small
modification of the above considerations gives the inequality
$\rho_E(r)\!\leq \!\frac{A}{c}\rho (1)\!+\!\log (Cr).$ The last
property in b) can be deduced from the assumption that $\rho$ is
doubly bounded and from the L-regularity of~$E$. \vskip 1mm

c) If $V_E(x+\zeta e_j)\leq \rho (|\zeta|)$ then for $1\leq \deg
P\leq n$ and $||P||_E=1$ we can write
$$\left\vert D^{ke_j}P(x)\right\vert\leq k!\inf\limits_{0<r\leq
1}r^{-k}\exp (n\sup\limits_{|\zeta |=r}V_E(x+\zeta e_j))$$
$$\leq k!\inf\limits_{0<r\leq
1}r^{-k}\exp (n\rho (r))=k!\exp (-k\psi^{-1}(k/n))\exp (n\rho
(\psi^{-1}(k/n))).$$ Hence
$$ \left\| D^{ke_j}P\right\|_E\leq k!\exp (-k\psi^{-1}(k/n))\exp (n\rho
(\psi^{-1}(k/n)))$$ and thus

$$ \left\| D^{\alpha }P\right\|_E\leq \alpha!\exp \left( -\sum\limits_{j=1}^{N} \ \!\!\! '\alpha_j\psi^{-1}
\left( \frac{\alpha_j}{n}\right)+n\sum\limits_{j=1}^{N} \ \!\!\!
'\rho\left( \exp\psi^{-1}\left(
\frac{\alpha_j}{n}\right)\right)\right)$$
$$\leq \alpha! \exp\left(-|\alpha| \psi^{-1}\left(c_N
|\alpha |/n\right) +n\rho\left(\exp\psi^{-1}\left( c_N|\alpha |/n
\right)\right)+ C_N|\alpha |\right).$$ We can see that the
assumptions of b) are satisfied and if we make use of it, then we
prove assertion  c). \hfill $\square$\vskip 3mm

As an application of Th.4.2.a we get the following bounds \vskip
3mm

 \noindent \textbf{Corollary 4.3.} a) If $\rho
(r)=Ar^\sigma,\ A\geq 1/\sigma$ then
$$||D^\alpha P||_E\leq \alpha!\left(A\sigma \sqrt{N}^{\sigma}e\right)^{|\alpha|/\sigma}
\!\!\left(
\frac{n}{|\alpha|}\right)^{|\alpha|/\sigma}\!\!\!||P||_E\leq\left(A\sigma
\sqrt{N}^{\sigma}e\right)^{|\alpha|/\sigma}\!\!\left(
\frac{1}{|\alpha|!}\right)^{1/\sigma -1}\!\!\!n^{|\alpha
|/\sigma}||P||_E.$$

b) If $\rho (r)=\frac{1}{s}(1/\log (e/r))^s$ then
$$||D^\alpha P||_E\leq \alpha ! ( \sqrt{N}/e)^{|\alpha |}\exp
\left((1+1/s)|\alpha |^{ \frac{s}{1+s}}n^{ \frac{1}{1+s}}\right)
||P||_E .$$

c) If $\rho (r)=r^\sigma\left( \log (1/r)+2/\sigma \right),\
\sigma\in (0,1],\ m=1/\sigma$ then
$$||D^\alpha P||_E\leq \alpha !N^{|\alpha |/2}\left(
\frac{e}{e-1}e^{2m}\right)^{|\alpha |} \left( \frac{n}{|\alpha
|}\right)^{m|\alpha|}\left(1+\log (n/|\alpha |)\right)^{m|\alpha|
}||P||_E.$$

\vskip 3mm

Th.4.2 generalizes implication $VMI\Rightarrow HCP$ in view of the
following \vskip 2mm

\noindent \textbf{Theorem 4.4.} {\it If $\sigma\in (0,1]$ and $
A\geq 1/\sigma$ then the function $\rho (r)=Ar^\sigma$ satisfies
conditions of a fit majorant and is {\it m}-, {\it mb}- and doubly
bounded.} \vskip 3mm

\noindent {\it Proof.} It is sufficient to check that $\rho$ is
{\it m}- and {\it mb}-bounded. Fix $c\in (0,1]$. Then $\psi
(t)=A\sigma e^{\sigma t}$, which implies $\psi^{-1}(ck/n)=\log
(ck/(A\sigma n))^{1/\sigma}$ and for $1\leq k\leq n$ we obtain
$$\exp\!\left(-k\psi^{-1}\!\!\left(
\frac{ck}n\right) \!+n\rho \!\left(\exp\psi^{-1}\!\!\left(
\frac{ck}n\right) \!\right)\right)=\left(
\frac{A\sigma}{c}\right)^{\!\!k/\sigma}\!\!\left(
\frac{n}{k}\right)^{\!\!k/\sigma}\!\!e^{ck/\sigma}\leq
\left(A\sigma e^c\frac{n}c\right)^{\!\!k/\sigma}\!\!\left(
\frac{1}{k!}\right)^{\!\!1/\sigma}\!\!\!\!.
$$ From this we conclude that
$$\log \left(1+\sum\limits_{k=1}^n\exp\left(-k\psi^{-1}\left(
ck/n\right) +n\rho \left(\exp\psi^{-1}\left(
ck/n\right)\right)\right)r^k\right)^{1/n}
$$
$$\leq \frac{1}{n}\log\left(1+\sum\limits_{k=1}^n\left( \left(A\sigma
e^cnr^\sigma/c\right)^k/k!\right)^{1/\sigma}\right)\leq
\frac{1}{n\sigma}\log\exp\left(A\sigma e^cnr^\sigma/c\right)
=\frac{A}{c}e^cr^\sigma\leq \frac{e}{c}\rho (r).
$$
It is easy to check that $\rho$ is {\it mb}-bounded  if we take
$c_s=1$ and $C_s=\frac{s}{\sigma}$. \hfill $\square$\vskip 3mm

The next theorem concerns another class of fit majorants.

\vskip 3mm

\noindent {\bf Theorem 4.5.} {\it The function $\rho
(r)=\frac{1}{p-1}\left(1/\log (e/r)\right)^{p-1}$ satisfies
conditions of a fit majorant and is {\it m}-, {\it mb}- and doubly
bounded.}

\vskip 3mm

To prove the above assertion we need the following crucial fact.
\vskip 3mm

 \noindent {\bf Lemma 4.6.} {\it Let $p,q>1$
and $ \frac{1}{p}+\frac{1}{q}=1$. There exists a positive constant
$B(p)$ such that for an arbitrary $s>0$ the following inequality
holds}
$$\sum\limits_{k=0}^\infty \exp (qk^{1/q})e^{-ks}\leq \exp
(B(p)s^{-(p-1)}).$$ \vskip 3mm

\noindent {\it Proof.} If $s\geq 1$ then
$$\sum\limits_{k=0}^\infty \exp (qk^{1/q})e^{-ts}\leq 1+e^{-s} \sum\limits_{k=1}^\infty
\exp (qk^{1/q})e^{-(k-1)}\leq 1+B_1(p)s^{-(p-1)}\leq \exp
(B_1(p)s^{-(p-1)}).$$ Now we will consider the more difficult case
$0<s<1$. We replace the series $\sum_{k=0}^\infty \exp (qk^{1/q})$
by the Laplace transform of the function $\exp (qt^{1/q})$.

We start with the simple observation that $\exp
(qk^{1/q})e^{-ks}\leq e^s\int\limits_k^{k+1}\exp
(qt^{1/q})e^{-ts}dt,$ which implies
$$ \sum\limits_{k=0}^\infty \exp (qk^{1/q})e^{-ks}\leq e
\int\limits_0^\infty\exp (qt^{1/q})e^{-ts}dt\leq \exp s^{-(p-1)}
\int\limits_0^\infty\exp (qt^{1/q})e^{-ts}dt.$$
$$ \int\limits_0^\infty\exp (qt^{1/q})e^{-ts}dt=q\int\limits_0^\infty
e^{qt}e^{-t^qs}t^{q-1}dt=qs^{-1}\int\limits_0^\infty \exp
(qt/s^{1/q})e^{-t^q}t^{q-1}dt.$$ Put
$b=s^{-\frac{1}{q}\frac{1}{q-1}}$. We can write
$\int\limits_0^\infty \exp
(qt/s^{1/q})e^{-t^q}t^{q-1}dt=I_1+\exp(qs^{-(p-1)})I_2$, where
$$I_1=\int\limits_0^b \exp
(qt/s^{1/q})e^{-t^q}t^{q-1}dt,\ I_2= \int\limits_0^\infty \exp
(qt/s^{1/q})e^{-(t+b)^q}(t+b)^{q-1}dt.$$ We have
$$I_1\leq \exp (qs^{-(p-1)})\int\limits_0^\infty
e^{-t^q}t^{q-1}dt=q^{-1}\exp (qs^{-(p-1)})<\exp (qs^{-(p-1)}).$$
Since $(b+t)^q=b^q+qt/s^{1/q}+\frac{1}{2}q(q-1)(b+\theta
t)^{q-2}t^2,\ \theta\in (0,1),$ we obtain
$$I_2=\exp (-s^{-(p-1)})\int\limits_0^\infty
\exp(-\frac{1}{2}q(q-1)(b+\theta t)^{q-2}t^2)(t+b)^{q-1}dt.$$ If
$1<q\leq 2$ then $(b+\theta t)^{q-2}\geq (b+t)^{q-2}$ and we get
$$\int\limits_0^\infty
\exp(-\frac{1}{2}q(q-1)(b+\theta t)^{q-2}t^2)(t+b)^{q-1}dt\leq
\int\limits_0^\infty \exp(-\frac{1}{2}q(q-1)(b+
t)^{q-2}t^2)(t+b)^{q-1}dt$$ $$=s^{-(p-1)}\int\limits_0^\infty
\exp(-\frac{1}{2}q(q-1)s^{-(p-1)}(t+ 1)^{q-2}t^2)(t+1)^{q-1}dt$$
$$ \leq s^{-(p-1)}\int\limits_0^\infty \exp(-\frac{1}{2}q(q-1)(t+
1)^{q-2}t^2)(t+1)^{q-1}dt=s^{-(p-1)}B_2(p)\leq\exp (
B_2(p)s^{-(p-1)}).
$$
If $q\geq 2$ then
$$\int\limits_0^\infty
\exp(-\frac{1}{2}q(q-1)(b+\theta t)^{q-2}t^2)(t+b)^{q-1}dt\leq
\int\limits_0^\infty
\exp(-\frac{1}{2}q(q-1)b^{q-2}t^2)(t+b)^{q-1}dt$$ $$\leq\!
(s^{-(p-1)})^{1/p}\!\!\int\limits_0^\infty
\exp(-\frac{1}{2}q(q-1)t^2)(t+1)^{q-1}dt\!\leq\!
(s^{-(p-1)})^{1/p}B_3(p)\!\leq\! \exp (s^{-(p-1)}B_3(p)^p/p).$$
The proof is completed by combining all particular cases. \hfill
$\square$ \vskip 3mm

\noindent {\it Proof of Theorem 4.5.} We have
$$1\!+\!\sum\limits_{k=1}^n\!\exp\!\!\left(\!\!-k\psi^{-1}\!\!\left( \frac{ck}n\right)
\!\!+n\rho \!\left(\exp\psi^{-1}\!\!\left(
\frac{ck}n\right)\!\right)\!\right)\!\!r^k\!\!=\!
1\!+\!\sum\limits_{k=1}^n\!\exp\!\! \left(\! \left(
\frac{1}{qc^{1/p}}\!+\!\frac{1}{p}c^{1/q}\!\right)\!qk^{1/q}n^{1/p}\!\right)\!\!\left(\frac{r}e\right)^{\!\!k}\!\!\!.$$
Observe that for $c\in (0,1]$ we have $
g(c):=\frac{1}{q}c^{-1/p}+\frac{1}{p}c^{1/q}\geq 1$. Applying
Lemma 4.6. with $s=u/(g(c)n^{1/p})$ we obtain
$$  1+\sum\limits_{k=1}^n\exp (qk^{1/q}n^{1/p}) \exp
(-ku/(g(c)n^{1/p})\leq \exp (B(p)g(c)^{p-1}n^{1/q}u^{-(p-1)}) .$$
Moreover, we see that
$$\left(1+\sum\limits_{k=1}^n\exp
(g(c)qk^{1/q}n^{1/p})e^{-ku}\right)^{1/g(c)n^{1/p}}\leq
1+\sum\limits_{k=1}^n\exp (qk^{1/q}n^{1/p}) \exp
(-ku/(g(c)n^{1/p}))$$
$$\leq \exp (B(p)g(c)^{p-1}n^{1/q}u^{-(p-1)})$$ which implies that
$$ 1+\sum\limits_{k=1}^n\exp
(g(c)qk^{1/q}n^{1/p})e^{-ku}\leq \exp (B(p)g(c)^pnu^{-(p-1)}).$$
If we now put $u=\log (e/r)$ for $r\in (0,1]$ then
$$\log \left( 1+\sum\limits_{k=1}^n\exp
(g(c)qk^{1/q}n^{1/p})(r/e)^k\right)^{1/n}\leq B(p)g(c)^p(\log
(e/r))^{-(p-1)}$$
$$=(p-1)B(p)( \frac{1}{q}+\frac{1}{p}c)^p \frac{1}{c}\frac{1}{p-1}(\log
(e/r))^{-(p-1)}\leq \frac{(p-1)B(p)}{c}\rho (r),$$ which gives
(\ref{m}).

We leave it to the reader to verify that  $\rho$ is {\it
mb}-bounded if we take $C_s=0$ and if $c_s\in (0,1]$ is chosen so
that
$$\text{ \ \ \ \ \ \ \ \ \ \ \ \ \ \ \ \
\ \ \ \ } q\sup\{ \sum\limits_{j=1}^s\lambda_j^{1/q}:\
\lambda_j\geq 0,\sum\limits_{j=1}^s\lambda_j=1\}\leq
c_s^{-1/p}(1+\frac{c_s}{p-1}). \text{ \ \ \ \ \ \ \ \ \ \ \ \ \ \
\ \  $\square$}$$ \vskip 3mm

We end this section with the following example.

\vskip 3mm

 \noindent \textbf{Example 4.7.} If $E$ is a unit ball in $
 \mathbb{C}^N$ then (cf. Section 2) $\rho_E(r)=\log (1+r/C(E))\leq r/C(E)$
 that is equivalent to $V_E(z)\leq \textrm{dist}(z,E)/C(E),\
 z\in\mathbb{C}^N$. We can take $\rho (r)=\max (1,1/C(E))r$ and thus  we get for $P\in \mathcal{P}(
 \mathbb{C}^N),\ \deg P\leq n,$
 $$||D^\alpha P||_E\leq (\max (1,1/C(E))
 (\sqrt{N}e)^{|\alpha|}n^{|\alpha|}||P||_E.$$

 \vskip 5mm

{\bf  5. HCP of compact subsets of $\mathbb{R}^N$}

\vskip 3mm

\noindent {\bf Remark 5.1.} For any set $E\subset\mathbb{R}^N$ it
is sufficient to consider only polynomials with real
coeffi\-cients. Indeed, if $P\in \mathcal{P}( \mathbb{C}^n)$ then
$P=Q+iR$ where $P,Q\in \mathcal{P}( \mathbb{R}^N),\ \deg P=\max
(\deg Q,\deg R)$ and
$$||P||_E=\sup\limits_{|\theta|\leq \pi}||\cos\theta
Q+\sin\theta R||_E,\ \ \ \ \ ||D^\alpha
P||_E=\sup\limits_{|\theta|\leq\pi}|| D^\alpha (\cos\theta
Q+\sin\theta R)||_E.$$ Hence, if we have the bound $||D^\alpha
P||_E\leq C(n,k)||P||_E$ for all $P\in \mathcal{P}(
\mathbb{R}^N),\deg P\leq n, |\alpha |\leq n$ then the same is true
for polynomials with complex coefficients.
\medskip

Observe that the following identity holds for real polynomials $P$
$$\Vert\text{\rm grad}\, P(x)\Vert_2^2=\frac{1}{2}\Delta
(P^2(x))-P(x)\Delta P(x)$$ and consequently,
$$E\in AMI(m,M)\ \Leftrightarrow\ \exists M'\ \forall
P\in\mathcal{P}( \mathbb{R}^N)\ ||\Delta P||_E\leq M'(\deg
P)^{2m}||P||_E.$$

Note also that, if $N=2$,  then
$$E\in AMI(m,M)\ \Leftrightarrow\ \exists M'\ \forall
P\in\mathcal{P}( \mathbb{R}^N) \ \ ||\frac{\partial}{\partial z}
P(x,y)||_E\leq M'(\deg P)^m ||P||_E$$
$$\Leftrightarrow\ \exists M'\ \forall
P\in\mathcal{P}( \mathbb{R}^N) \ \
||\frac{\partial}{\partial\overline{z}} P(x,y)||_E\leq M'(\deg
P)^m ||P||_E$$ where, as usual,  $$\frac{\partial}{\partial z}
P(x,y)=\frac{1}{2} \left( \frac{\partial P(x,y)}{\partial
x}-i\frac{\partial P(x,y)}{\partial y}\right),\
\frac{\partial}{\partial \overline{z}} P(x,y)=\frac{1}{2} \left(
\frac{\partial P(x,y)}{\partial x}+i\frac{\partial
P(x,y)}{\partial y}\right) .$$

For compact subsets of $ \mathbb{R}^N$ we can take only real
polynomials in the definition of Siciak's extremal function (cf.
\cite{[1]}):
$$\Phi_E(z)=\sup\left\{ |P(z)|^{1/\deg P}:\ P\in \mathcal{P}(
\mathbb{R}^N),\ \deg P\geq 1,||P||_E\leq 1\right\},\
z\in\mathbb{C}^N .$$ \vskip 3mm

The following result is a consequence of \cite[Th.2.4]{[3]}

\vskip 3mm

\noindent \textbf{Proposition 5.2.} \textit{If $E$ is a compact
subset of $\mathbb{R}^N$ then \be V_E(x+iy)\leq
\frac{1}{\pi}\int\limits_{-\infty}^{+\infty}V_E(x+ty)
\frac{dt}{1+t^2} \label{wst} \ee for every
$z=x+iy\in\mathbb{C}^N$. Equality holds in {\rm (\ref{wst})} if
$N=1$ (for any $z\in\mathbb{C}$).}

\vskip 3mm

As a corollary of Prop.5.2 we obtain

\vskip 3mm

\noindent \textbf{Theorem 5.3.} \textit{Let $E$ be a compact set
in $\mathbb{R}^N$. Assume that for every $x\in \mathbb{R}^N$ the
inequality holds \be V_E(x)\leq B(\textrm{dist}(x,E))^\gamma
\label{zal} \ee where $B>0,\ \gamma\in
 (0,1)$ are constants independent of $x$. Then for all
 $z\in\mathbb{C}^N$
 $$V_E(z)\leq \widetilde{B}(\textrm{dist}(z,E))^\gamma,\ \ \ {\it
 with} \ \ \
 \widetilde{B}=\frac{B}{\pi}\int\limits_{-\infty}^{+\infty}(1+t^2)^{\gamma /2-1}dt=
 \frac{B}{\sqrt{\pi}} \frac{\Gamma (1/2-\gamma /2)}{\Gamma (1-\gamma
 /2)}.$$}

\vskip 3mm

\noindent \textit{Proof.}  Evidently, for $z=x+iy\in \mathbb{C}^N$
we have
 $$ \textrm{dist}(z,E)=( \textrm{dist}(x,E)^2+\| y\|_2^2)^{1/2}.$$
 Inequality (\ref{zal}) is equivalent to
 $$V_E(x)\leq B\Vert x-x_0\Vert_2^\gamma  \ \ \ \ {\rm for \ all }
 \ \ \ x_0\in E.$$ By Prop.5.2, we get
$$V_E(z)=V_E(x+iy)\leq
\frac{1}{\pi}\int\limits_{-\infty}^{+\infty}V_E(x+ty)
\frac{dt}{1+t^2}\leq
\frac{B}{\pi}\int\limits_{-\infty}^{+\infty}\Vert
x-x_0+ty\|_2^\gamma \frac{dt}{1+t^2}$$
$$\leq \frac{B}{\pi}\int\limits_{-\infty}^{+\infty}(\Vert
x-x_0\|_2+|t|\|y\|_2)^\gamma \frac{dt}{1+t^2}\leq
\frac{B}{\pi}\int\limits_{-\infty}^{+\infty}(\Vert
x-x_0\|_2^2+\|y\|_2^2)^{\gamma /2} (1+t^2)^{\gamma /2}
\frac{dt}{1+t^2}$$
$$= \widetilde{B}(\Vert
x-x_0\|_2^2+\|y\|_2^2)^{\gamma /2}.$$ As $x_0$ is an arbitrary
point of $E$, we obtain $V_E(z)\leq
\widetilde{B}(\textrm{dist}(z,E))^\gamma$, which completes the
proof. \hfill $\square$

\vskip 3mm

It may be worth reminding the reader that if a compact set
$E\subset \mathbb{R}^N$ admits the A.Markov inequality then the
exponent $m$ in (\ref{AMI}) is at least equal to 2 (see e.g.
\cite{[11]}). Therefore, the exponent $\gamma$ in (\ref{zal}) may
be at most equal to $\frac12$.

\vskip 3mm

\noindent \textbf{Corollary 5.4.} \textit{If for every $x\in
\mathbb{R}^N$ the inequality holds \be V_E(x)\leq
B(\textrm{dist}(x,E))^\gamma \label{zal2} \ee with $B>0,\
\gamma\in
 (0,\frac12]$ independent of $x$, then for all
$z\in\mathbb{C}^N$
$$ V_E(z)\leq \widetilde{B}(\textrm{dist}(z,E))^\gamma , \ \ \ \ \widetilde{B}=
\frac{B}{\sqrt{\pi}} \frac{\Gamma (1/4)}{\Gamma (3/4)}.$$}

\vskip 3mm

\noindent \textbf{Corollary 5.5.} \textit{If $E$ is a compact
subset of $\mathbb{R}^N$ and $\gamma\in (0,1]$, then the following
conditions are equivalent: }

\vskip 2mm

$(i) \ \ \ E\in HCP(\gamma, B_1) \ \ \ {\it with \ some} \ \ \
B_1\ge 1,$

\vskip 2mm

$(ii) \ \, \ {\it inequality \ (\ref{zal}) \ holds \ for \ all} \
\ x\in \mathbb{R}^N \ \ {\it with \ some} \ \ B_2\ge 1$ {\it
independent of $x$,}

\vskip 2mm

$(iii) \ E\in VMI(\frac1\gamma,B_3) \ \ \ {\it with \ some} \ \ \
B_3\ge 1,$

\vskip 2mm

$(iv)$ \  \ {\it inequality (\ref{VMI}) holds with some \ $B\ge 1$
for all polynomials $P$ of real coefficients. }

\vskip 5mm

\noindent \textbf{Example 5.6.} For $E=[-1,1]$ we have
$V_E(x)=\log h(\max (1,|x|))$ where $h(t)=t+\sqrt{t^2-1}$ if
$t\geq 1$. If $x_0\in E$ then
$$V_E(x)\leq \log h(1+|x-x_0|)=\log \left( 1+|x-x_0|^{1/2}(|x-x_0|^{1/2}+(|x-x_0|+2)^{1/2})\right).$$
Since \ $\log (1+t)\leq \frac{1}{\alpha}t^{\alpha}$ \ for $t\geq
0$, $0<\alpha \leq 1$, we obtain \[ V_E(x)\leq \left\{
       \begin{array}{cc}
        (1+\sqrt{3}) |x-x_0|^{1/2} \ \ \ \ \ \ & ; \ |x-x_0|\leq
1,\\
        2(1+\sqrt{3})^{1/2}|x-x_0|^{1/2} & ; \ |x-x_0|>
1,
       \end{array}
       \right. \ \ \ \ \ \ \ \ \ \ \ \ \ \ \ \ \ \ \ \ \ \ \ \ \ \ \ \]
hence $V_E(x)\leq 2(1+\sqrt{3})^{1/2}|x-x_0|^{1/2}$ for all $x$,
$x_0$ and thus
$$V_E(x)\leq 2(1+\sqrt{3})^{1/2}(\textrm{dist}(x,E))^{1/2}.$$

\vskip 3mm

\noindent \textbf{Example 5.7.} Let $E$ be a convex body in
$\mathbb{R}^N$ that is not symmetric with respect to the origin.
Fix $\xi\in S^{N-1}$ and put
$$a_\xi (E)=\min\limits_{x\in E}\langle x,\xi\rangle ,\ \ \
b_\xi (E)=\max\limits_{x\in E}\langle x,\xi\rangle ,\ \ \ \rho_\xi
(E)=b_\xi (E)-a_\xi(E).$$ The last value is called the width of
$E$ in the direction $\xi$. The minimal width of $E$ is given by
$\omega (E)=\inf\limits_{\xi\in S^{N-1}}\rho_\xi (E)$. For
$x\in\mathbb{R}^N$ it follows that (see \cite{[9]})
$$V_E(x)=\sup\limits_{\xi\in
S^{N-1}}V_{[a_\xi (E),b_\xi (E)]}(\langle x,\xi\rangle )
=\sup\limits_{\xi\in S^{N-1}}V_{[-1,1]}(2\langle
x,\xi\rangle/\rho_\xi (E)-(b_\xi (E)+a_\xi (E))/\rho_\xi (E))$$
Therefore, in the same way as in Example 5.6, we have
$$V_E(x)\leq \sup\limits_{\xi\in
S^{N-1}}\log h\left( 1+2|\langle x-x_0,\xi\rangle |/\rho_\xi
(E)\right)\leq \log h(1+2\| x-x_0\|_2/\omega (E))$$ $$\leq
2(1+\sqrt{3})^{1/2}(2\| x-x_0\|_2/\omega (E))^{1/2}$$ for any
fixed $x_0\in E$ and in consequence we get
$$V_E(x)\leq
(2+2\sqrt{3})^{1/2}(4\textrm{dist}(x,E)/\omega(E))^{1/2}\leq
(2+2\sqrt{3})^{1/2}(\textrm{dist}(x,E)/C(E))^{1/2},$$ where $C(E)$
is the L-capacity of $E$ (see \cite[Example 3.4]{[4]}). In
particular, there exists an absolute constant $\mathcal{B}$ such
that for all dimensions $N$ and for all convex bodies
$E\subset\mathbb{R}^N$ the inequality holds
$$V_E(z)\leq  \mathcal{B}(\textrm{dist}(z,E)/C(E))^{1/2},\ z\in\mathbb{C}^N.$$
By Th.2.12, we can deduce that these sets belong to $VMI(2,
\sqrt{N}\mathcal{B}^2e^2/[4C(E)])$.

\vskip 3mm

Now we recall a definition of a class of UPC sets introduced by
Paw\l{}ucki and Ple\'sniak \cite{[16]} who have shown its
importance in approximation theory. In particular, they have
proved a deep result (cf. \cite[Cor. 6.5]{[16]}) that every fat
compact subanalytic subset of $ \mathbb{R}^N$ belongs to this
class (see also \cite{pi1}).

\vskip 1mm

 Let $s\geq 1$, $S>0$ and $d\in \{1,2,\ldots\}$.

\vskip 2mm

\noindent \textbf{Definition 5.8.}  A compact set $E\subset
\mathbb{R}^N$ is called \textit{uniformly polynomially cuspidal}
 ($E\in UPC(s,S,d)$ in short) if for every $x_0\in E$ we can find a polynomial mapping
$\varphi \, : \, \mathbb{R} \rightarrow \mathbb{R}^N$ of degree at
most $d$ such that $\varphi (1)=x_0$ and $$\textrm{dist}(\varphi
(t),\mathbb{R}^N\setminus E)\geq S(1-t)^s \ \ \ \ {\rm for} \ \ \
\ t\in [0,1].$$

\vskip 3mm

It is rather difficult to find the optimal constant $s$ in the
last inequality. However, calculations are much simpler for the
following modification of the above definition.

\vskip 3mm

\noindent {\bf Definition 5.9.}(cf. \cite{[2]}) Let  $v$ be a
fixed unit vector in $ \mathbb{R}^N$. A compact set $E\subset
\mathbb{R}^N$ is called \textit{uniformly polynomially cuspidal in
direction $v$}
 ($E\in UPC_v(s,S,d)$ in short) if for every $x_0\in E$ we can find a polynomial mapping
$\varphi \, : \, \mathbb{R} \rightarrow \mathbb{R}^N$ of degree at
most $d$ such that $\varphi (1)=x_0$ and $$\textrm{dist}_v(\varphi
(t),\mathbb{R}^N\setminus E)\geq S(1-t)^s \ \ \ \ {\rm for} \ \ \
\ t\in [0,1].$$ Here $\textrm{dist}_v(x,\mathbb{R}^N\setminus
E):=\sup\{r\geq 0:\ [x-rv,x+rv]\subset E\}$. \vskip 3mm

If $E\in UPC(s,S,d)$ then $ E\in UPC_v(s,S,d)$ for every unit
vector $v$. An open problem is whether conditions $ E\in
UPC_{v_j}(s_j,S_j,d_j),\ j=1,\dots ,N$,  $v_1,\dots ,v_N$ that are
linearly independent imply $E\in UPC(s,S,d)$ with some $S,s,d$.
 It seems that this may not be true for $N\geq 3$. However, as
 an application of the
proposition given below, we prove that it $E\in
UPC_{e_j}(s_j,S_j,d),\ j=1,\dots ,N$, where $(e_j)_j$ is the
canonical basis,  then we get $E\in HCP(\frac1{2s})$, where
$s=\max\limits_{1\leq s\leq N}s_j$. In particular, if $E\in
UPC(s,S,d)$ then $E\in HCP( \frac{1}{2s})$ that essentially
improves earlier result by Paw\l{}ucki and Ple\'sniak
\cite[Th.4.1]{[16]} (see also \cite{[17]}). As a corollary we get
a wide class of non-convex sets that satisfy $VMI$.
\vskip 3mm

\noindent {\bf Theorem 5.10.} {\it If $E\in UPC_v(s,S,d)$ and
$\varepsilon_0\in (0,1)$ then there exists
$C_0=C_0(\varepsilon_0)>0$ such that for every $|\zeta|\leq r_0=
\frac{S}{\sqrt{2}}(1-\varepsilon_0)^s$ the inequality holds
$$V_E(x+\zeta v)\leq C_0|\zeta |^{1/(2s)}.$$}

\noindent \textit{Proof.} Let $L_0=\sqrt{2}/S$. Put $g(\zeta
)=\frac{1}{2}(\zeta +\zeta^{-1}),\ \widehat{g}(\zeta
)=\frac{1}{2}(\zeta -\zeta^{-1})$. For $\rho > 1$ write $a=a(\rho
)=((1+g(\rho ))/2)^{-1},\ b=b(\rho )=((g(\rho )-1)/2)^{-1},\
c=c(\rho )=1/\widehat{g}(\rho )$. We have $b=a(1-a)^{-1},\
c=\frac{1}{2}a(1-a)^{-1/2}$. Fix
$\zeta_0=\alpha_0+i\beta_0\in\mathbb{C}^N$ and  $x_0$ in $E$ and
put $u_0=\alpha_0v,\ v_0=\beta_0v$. Define
$$\psi (\zeta )=\varphi (a(\rho )\frac{1}{2}(g(\zeta )+1))+\frac{1}{2}(g(\zeta
)-1)b(\rho )u_0+i\widehat{g}(\zeta )c(\rho )v_0, \ \ \ \
\zeta\in\mathbb{C}, \ \ |\zeta|\geq 1$$ where $\varphi$ is a
polynomial mapping chosen to $x_0$ by the definition of the $UPC$
property.

Since
$$\psi (e^{i\theta})=\varphi (a(\rho)r\tau )+(\tau -1)b(\rho )u_0\pm
2\sqrt{\tau (1-\tau )}c(\rho )v_0 \ \ {\rm for} \: \tau\!=\!
\frac{\cos\theta+1}{2}, \: \theta\!\in\![0,2\pi ],$$ we have $\psi
(e^{i\theta})\in E$ whenever
$$ (1-\tau )b(\rho )|\alpha_0| +2\sqrt{\tau (1-\tau)}c(\rho )|\beta_0|\leq
S(1-ar)^s,\ \tau \in[0,1].$$ The last condition is satisfied if
$$a\left( \frac{1-\tau }{1-a}|\alpha_0| +\sqrt{\tau\frac{1-\tau}{1-a}}|
\beta_0| \right)\leq S(1-a\tau )^s$$ and consequently, if
$$a\left( \frac{1-\tau}{1-a\tau}(1-a\tau )^{-(s-1)}\frac{|\alpha_0|}{1-a}
+\sqrt{\frac{\tau (1-\tau )}{1-a\tau}}(1-a\tau )^{-(s-1/2)}
\frac{|\beta_0|}{\sqrt{1-a}}\right)\leq S.$$ Since $ a,\tau \in
[0,1]$,
 $1-\tau \leq 1-a\tau$ and  $ |\alpha_0|+|\beta_0|\leq
 \sqrt{2}|\zeta_0|$, the last condition
 will hold if
 $\sqrt{2}|\zeta_0| \leq S(1-a)^s.$

 Assuming $L_0|\zeta_0|\leq (1-\varepsilon_0)^s$ and taking $ \rho =
  h\left( \frac{1+(L_0\delta )^{1/s} }{1-(L_0\delta
  )^{1/s}}\right),\ h(t)=t+\sqrt{t^2-1}$ we have $L_0|\zeta_0| =(1-a(\rho
  ))^s$, that is $\psi (\{|z|=1\})\subset E$.
By the maximum principle for subharmonic functions (applied to the
domain $\{z\in \mathbb{C}\, : \, |z|>1\}$), we get
  $$ \log V_E(\psi (\zeta ))\leq d\log |\zeta |,\ |\zeta |\geq
  1.$$ In particular,
  $$V_E(x_0+\zeta_0v)=V_E(\psi (\rho))\leq d\log\rho =
  d\log h\left( \frac{1+(L_0|\zeta_0| )^{1/s} }{1-(L_0|\zeta_0|
  )^{1/s}}\right).$$

The inequality $1-(L_0|\zeta_0|
  )^{1/s}\geq \varepsilon_0$ implies that
  $$h\left( \frac{1+(L_0\delta )^{1/s} }{1-(L_0\delta
  )^{1/s}}\right)\leq h(1+(2/\varepsilon_0)(L_0\delta )^{1/s})
   \leq 1+A\delta^{1/(2s)}$$
  where $$ A=(2/\varepsilon_0)^{1/2}L_0^{1/(2s)}(
  \sqrt{2/\varepsilon_0}+\sqrt{2(1-\varepsilon_0)/\varepsilon_0})
  \leq (4/\varepsilon_0) L_0^{1/(2s)}=B.$$ Since for every  $d\geq 1$
  the function $((1+x)^d-1)/x$  is increasing for $x>0$, we obtain
  $$(1+B|\zeta_0| ^{1/(2s)})^d\leq 1+C_0|\zeta_0|^{1/(2s)}$$
  where  $C_0=\max\limits_{r\leq
  r_0}\left((1+Br^{1/(2s)})^d-1\right)/r^{1/(2s)}=
  \left((1+Br_0^{1/(2s)})^d-1\right)/r_0^{1/(2s)}$.

  Finally $V_E(x_0+\zeta_0v)\leq \log
  (1+C_0|\zeta_0|^{1/(2s)})\leq
  C_0|\zeta_0|^{1/(2s)}$. \hfill $\Box$
\vskip 3mm

 Applying Cor.2.13 and Remark 3.2 we get the
following result which specifies an earlier result by Paw\l{}ucki
and Ple\'sniak (cf. \cite[Th.2.1]{[16]}). \vskip 3mm

 \noindent {\bf Corollary 5.11.} {\it  If $E\in
 UPC_{e_j}(s_j,S_j,d_j),\ j=1,\dots ,N$ then there exists a
 constant $B$ such that $E\in HCP(\gamma ,B)$ with $\gamma
 =1/(2\min\limits_js_j)$. In particular, if $E\in UPC(s,S,d)$ then $E\in
 HCP(1/(2s),B)$. }

\vskip 5mm

{\bf  6. Applications of Theorem 2.12 for disconnected sets.}

\vskip 3mm

The first proposition regards certain onion type sets in the
complex plane that are very useful in a problem concerning local
and global Markov's properties (see L.Bialas-Ciez and R.Eggink,
{\it Equivalence of the global and local Markov inequalities in
the complex plane}, in preparation).

\vskip 3mm

\noindent \textbf{Proposition 6.1.} \textit{Let $(a_j)_j$ be a
strictly decreasing sequence of positive numbers such that $a_1 =1
$, $a_j\rightarrow 0$ as $j\rightarrow \infty$ and let
$\varphi_j\in (0,\frac\pi2)$ for $j=1,2,\ldots$. Put
\[ C_j:=\{ a_je^{it} \: : \: t\in [\varphi_j,2\pi]\} \ \ \ \ \ \ {
for} \ \ j=1,2,\ldots,\]
\[ E:= \{0\} \cup \bigcup_{j=1}^\infty \: C_j.\]
If \ $|1-e^{i\varphi_j}|\le a_{j+1}$ \ for \ $j=1,2,...$ \ then \
$E\in HCP(\frac16,B)$ for some $B>0$.}

\vskip 3mm

\noindent {\it Proof.} First, we note that
\[ F \ := \ \{ e^{i\,t} \: : \: t\in [\pi/2, 2\pi]\}\]
is a connected compact set and so $F\in HCP(\frac12,B_F)$ with
some constant $B_F\ge 1$ (see e.g. \cite[Cor.2.2]{[8]}). From
implication (\ref{hcpvmi}) we see that $F\in VMI(2,M_F)$ with
$M_F=(e{B_F}/2)^2$. We can assume that $M_F\ge \max\{2e,
1/C(E)\}$.

For any polynomial $P$ of degree at most $n$, for $k\in
\{1,\ldots,n\}$ and $z_0\in E$, we will prove the inequality \be
|P^{(k)}(z_0)| \ \le \ M^k \: \frac{n^{6k}}{k^{5k}} \: \|P\|_E \ \
\ \ \ {\rm where} \ \ \ M= 3M_F
\exp\left({3M_F(1+e^{3M_F})}\right). \label{wyst1}\ee By Th.2.12
and since $k^k\ge k!$, condition (\ref{wyst1}) implies that $E\in
HCP(\frac16,B)$ with $B=6\, M^{1/6}$. Therefore, the proof is
completed by showing (\ref{wyst1}).

Observe first that for any monic polynomial $P$ of degree $n$ and
for $k=n$ we have
\[ \|P^{(k)}\|_E \ = \ \left(\frac1{C(E)}\right)^k  n! \
(C(E))^n \ \le \ \left(\frac1{C(E)}\right)^k  n! \ \|P\|_E,\]
because $C(E)$ is equal to the Chebyshev constant of $E$.
Consequently,
\[ \|P^{(k)}\|_E \ \le \ \left(M_F \right)^k \frac{n^{6k}}{k^{5k}}
 \ \|P\|_E\] for all polynomials of degree at most $n$ not necessary monic.
Thus condition (\ref{wyst1}) is fulfilled for $k=n$.

Consider now $k<n$. We first examine $z_0=0$ and we will show that
\be |P^{(k)}(0)| \ \le \ \left( \frac{3M_F\, e^{3M_F}\, n^4}{k^3}
\right)^k  \, ||P||_E. \label{1}\ee  For this purpose, find $j\in
\mathbb{N}$ such that $\frac1{a_j}\! \le \!\frac{n^2}{k^2} \!<
\!\frac1{a_{j+1}}$. By Cauchy's integral formula, \[ |P^{(k)}(0)|
\ \le \ \frac{1}{a_j} \: ||P^{(k-1)}||_{C(0,a_j)} \ ,\] where
$C(0,a_j)$ is the circle with the radius $a_j$ about the origin.
The norm of $P$ on ${C(0,a_j)}$ is attained at some point $w_0\in
C(0,a_j)$. Put $F_j:=a_j F$. Obviously, $F_j\in VMI(2, M_F/a_j)$.
If $w_0\in F_j$, we have
\[ |P^{(k)}(0)| \le \frac{1}{a_j} \ ||P^{(k-1)}||_{F_j}
\le  \left( \frac{1}{a_j}\right)^{\!\!\!k} \! \frac{\left(M_F\:
n^2\right)^{k-1}}{(k-1)!} \ ||P||_{F_j} \le  \left(
\frac{1}{a_j}\right)^k \left( \frac{3M_F\: n^2}{k-1}\right)^{k-1}
\: ||P||_{F_j} \] \[ \le \ \left( \frac{n^2}{k^2}\right)^k \left(
\frac{3M_F\: n^2}{k}\right)^{k} \: ||P||_{F_j} \ \le \  \left(
\frac{3M_F\: n^4}{k^3}\right)^{k} \: ||P||_{E} \ \ \ \ \ \ \ {\rm
for} \ \ k\ge 2,
\] because $t\mapsto \left(
\frac{M_F\, n^2}{t}\right)^{t}$ is an increasing function whenever
$t\in (0,n+1] \subset (0,\frac{M_F\, n^2}{e}]$ and $F_j\subset
C_j\subset E$. The case of $k=1$ is easy to verify and thus
inequality (\ref{1}) is fulfilled for all $k\in \{1,\ldots,n\}$,
$w_0\in F_j$.

If $w_0\in C(0,a_j)\setminus F_j$, by Taylor's formula and VMI for
$F$, we get
\[ |P^{(k)}(0)| \ \le \ \frac{1}{a_j} \ |P^{(k-1)}(w_0)| \ \le \ \frac{1}{a_j}
\sum_{l=0}^{n-k+1} \frac1{l!} \,|P^{(k-1+l)}(a_j)| \,
|a_j-w_0|^l\] \[ \le \frac{1}{a_j} \sum_{l=0}^{n-k+1} \frac1{l!}
\, \left( \frac{1}{a_j}\right)^{k-1+l} \left( \frac{3M_F\:
n^2}{k-1+l}\right)^{k-1+l} \: ||P||_{F_j} \ a_j^l|1-e^{i\,
\varphi_j}|^l \]
\[ \le \left( \frac{1}{a_j}\right)^{k} \ \sum_{l=0}^{n-k+1}
\frac1{l!} \,  \left( \frac{3M_F\: n^2}{k+l}\right)^{k+l}
 a_{j+1}^l \: ||P||_{F_j} ,\] the last inequality being a consequence
of the assumption of Prop.6.1. Since $F_j\subset E$ and
$\frac1{a_j}\! \le \!\frac{n^2}{k^2} \!< \!\frac1{a_{j+1}}$, we
can write
\[ |P^{(k)}(0)| \le \left( \frac{n^2}{k^2}\right)^{\!\!k} \!\left(
\frac{3M_F\: n^2}{k}\right)^{\!\!k} \: \sum_{l=0}^{n-k+1}
\frac1{l!} \!  \left( \frac{3M_F\: n^2}{k}\right)^{\!\!l} \!\left(
\frac{k^2}{n^2}\right)^{\!\!l} ||P||_{E} \le  \left( \frac{3M_F \,
e^{3M_F} \: n^4}{k^3}\right)^{\!\!k} \! ||P||_{E}\] and this
yields inequality (\ref{1}).

We now turn to the case \ $z_0\ne 0$. Clearly, $z_0\in C_j$ \ for
some $j\in \{1,2,\ldots,\}$. If \ $a_j \le \frac{k^4}{n^4}$ then
by (\ref{1}) we have
\[ |P^{(k)}(z_0)| \le \sum_{l=0}^{n-k} \frac1{l!} \,|P^{(k+l)}(0)| \,
|z_0|^{l} \le \sum_{l=1}^n \frac1{l!} \, \left( \frac{3M_F \,
e^{3M_F} \: n^4}{(k+l)^3}\right)^{k+l} \, ||P||_{E} \: a_j^{l}
\] \[ \le \left( \frac{3M_F \,
e^{3M_F} \: n^4}{k^3}\right)^{k} \sum_{l=1}^n \frac1{l!} \, \left(
\frac{3M_F \, e^{3M_F} \: n^4}{k^3}\right)^{l}  \:
\left(\frac{k^4}{n^4}\right)^{l} \, ||P||_{E} \] \[ \le  \left(
\frac{3M_F \, e^{3M_F} \, \exp({3M_F \, e^{3M_F}}) \:
n^4}{k^3}\right)^{k}
 \,
||P||_{E} \] and (\ref{wyst1}) is proved in this case.

It remains to show estimate (\ref{wyst1}) if \ $a_j >
\frac{k^4}{n^4}$. Let $F_j'$ be a set obtained by a rotation of
$F_j$ about the origin such that \ $z_0\in F_j'\subset C_j$. Since
$F_j'\in VMI(2, {M_F}/{a_j})$, we have
\[ |P^{(k)}(z_0)| \le \|P^{(k)}\|_{F_j'} \le \frac1{a_j^k} \left(
\frac{3M_F \: n^2}{k}\right)^{k} \|P\|_{F_j'} \le \left(
\frac{n^4}{k^4}\right)^k \left( \frac{3M_F \: n^2}{k}\right)^{k}
\|P\|_{E}
\] and (\ref{wyst1}) is proved at every point $z_0\in
E$. \hfill $\Box$

\vskip 3mm

The second example presents the application of Th.2.12 for certain
compact sets consisting of infinitely many pairwise disjoint
subsets of $\mathbb{C}^N$.

\vskip 2mm

\noindent \textbf{Proposition 6.2.} \textit{Let $\mu\ge 2$, $b\in
(0,\sqrt{2}-1)$ and let $(a_j)_j$, $(r_j)_j$ be sequences of
positive numbers such that \[ a_1=2, \ \ r_1=1, \ \ \
a_j=r_j+r_j^2, \ \ \ r_j=b\,r_{j-1}^\mu \  for \ \ j\ge 2.\] Then
the set $E$ defined by \[ E := \{0\} \cup \bigcup_{j=1}^\infty \:
E_j, \ \ \ E_j:=\{ z=(z_1,\ldots,z_N)\in \mathbb{C}^N \: : \:
|z_1-a_j|\le r_j, |z_2|\le r_j, \ldots, |z_N|\le r_j\}.\] admits
the H\"older continuity property of the pluricomplex Green
function $HCP(\frac1{2\!+\!\mu},\! B)$ with some $B>0$.}

\vskip 1mm

\noindent {\it Proof.} Fix $n\in \{1,2,\ldots\}$ and a polynomial
$P$ of degree at most $n$. As a first step we shall show that for
each $\alpha\in \mathbb{N}_0^N$, $|\alpha|\le n$  \be
|D^{(\alpha)}P(0)| \ \le \ \left( \frac{e^{N}}{b}
\right)^{|\alpha|} \, \left( \frac{n^{1+\mu}}{|\alpha|^\mu}
\right)^{|\alpha|} \, ||P||_E. \label{2}\ee For this purpose, find
$j\ge 2$ such that $r_j\! < \!\frac{|\alpha|}{n} \!\le \!r_{j-1}$
where $|\alpha|\le n$ is fixed. From (\ref{3}) we have
\[ \|D^\alpha P\|_{E_j} \ \le \ \frac{n^{|\alpha|}}{r_j^{|\alpha|}} \ \|P\|_{E_j}\]
and thus, by Th.2.12 and Example 2.9, $E_j\in HCP(1,
\frac{N}{r_j})$. In particular, we get
\[ V_{E_j} (0) \ \le \ \frac{N}{r_j} \: {\rm dist}\, (0,E_j) \ = \ \frac{N}{r_j}
\: (a_j -r_j) \ = \ N \, r_j \ < \ \frac{|\alpha| \, N}{n}.\]
Formula (\ref{V_E}) leads us to
\[ |D^{(\alpha)}P(0)| \ \le \
\left( e^{V_{E_j}(0)} \right)^{n-|\alpha|} \,  ||D^\alpha
P||_{E_j}.\] By the above, it follows that
\[ |D^{(\alpha)}P(0)| \ \le \
e^{|\alpha|\, N}  \  \frac{n^{|\alpha|}}{r_j^{|\alpha|}} \ \:
\|P\|_{E_j} \ \le \ e^{|\alpha|\, N}  \ \frac{n^{|\alpha| + \mu
|\alpha|}}{(b|\alpha|^\mu)^{|\alpha|}} \ \: \|P\|_{E} \] and
inequality (\ref{2}) is proved.

Now consider $z_0\in E\setminus \{0\}$ and $|\alpha|\le n$. If
$z_0\in E_j$ and $r_j \ge \left(\frac{|\alpha|}{n}
\right)^{1+\mu}$ then \be |D^{(\alpha)}P(z_0)| \ \le \
\frac{n^{|\alpha|}}{r_j^{|\alpha|}} \ \: \|P\|_{E_j} \ \le \
\frac{n^{2|\alpha| + \mu |\alpha|}}{|\alpha|^{(1+\mu)|\alpha|}} \
\: \|P\|_{E} .\label{5} \ee In the case of $r_j <
\left(\frac{|\alpha|}{n} \right)^{1+\mu}$, Taylor's formula and
inequality (\ref{2}) yield
\[ |D^{(\alpha)}P(z_0)| \ \le \ \sum_{|\beta|\le n-|\alpha|}
\frac1{\beta!} \ |D^{\alpha+\beta} P(0)| \ \:
\|z_0\|_2^{|\beta|}\] \[ \le \ \sum_{|\beta|\le n-|\alpha|}
\frac1{\beta!} \left( \frac{e^{N}}{b} \right)^{|\alpha|+|\beta|}
\, \left( \frac{n^{1+\mu}}{(|\alpha|+|\beta|)^\mu}
\right)^{|\alpha|+|\beta|} \, ||P||_E \: \
(r_j\sqrt{N+8})^{|\beta|}\] \[ \le \ \left( \frac{e^{N}}{b}
\right)^{|\alpha|} \, \left( \frac{n^{1+\mu}}{|\alpha|^\mu}
\right)^{|\alpha|} \, ||P||_E   \sum_{l=0}^{n-|\alpha|}
\frac{N^l}{l!} \left( \frac{e^{N}\sqrt{N+8}}{b} \right)^{l} \,
\left( \frac{n^{1+\mu}}{|\alpha|^\mu} \right)^{l} \
\left(\frac{|\alpha|}{n} \right)^{(1+\mu)l}\] \[ \le \ \left(
\frac{1}{b} \ e^{N+N\sqrt{N+8}\:e^N/b}\right)^{\!\!|\alpha|} \!
\left( \frac{n^{1+\mu}}{|\alpha|^\mu} \right)^{\!\!|\alpha|} \,
||P||_E \ \le \ \left( \frac{1}{b} \
e^{N+N\sqrt{N+8}\:e^N/b}\right)^{\!\!|\alpha|} \!
\frac{\left(n^{1+\mu}\right)^{|\alpha|}}{(|\alpha|!)^\mu}  \,
||P||_E.
\] Hence and from inequalities (\ref{2}, \ref{5}) we conclude that
$E\in VMI(2+\mu, \frac{1}{b} \ e^{N+N\sqrt{N+8}\:e^N/b})$, and
Th.2.12 leads to $E\in HCP(\frac1{2+\mu},\! B)$ with
$B=\left(\frac{N}{b} \
e^{N+N\sqrt{N+8}\:e^N/b}\right)^\frac1{2+\mu} (2+\mu)$, which
proves the assertion. \hfill $\Box$

\vskip 4mm

\noindent \textbf{Remark 6.3.} We close this paper by offering two
questions for further research:

\noindent 1. Does the continuity of the pluricomplex Green's
function $V_E$ with respect to each variable separately imply the
L-regularity of $E$?

\noindent 2. Has the pluricomplex Green's function $V_E$ of a
Markov set $E$ the continuity property with respect to each
variable separately?

\noindent For the univariate case the answer to the second
question is partially known because if $E\subset \mathbb{R}$ then
it is L-regular (see \cite{[7]}).

\end{document}